\documentclass[10pt,twocolumn,twoside]{IEEEtran}
\usepackage{amsmath,amssymb,amsthm,epsfig,dsfont,color,subfigure,empheq,graphicx}
\usepackage{enumerate,url,algpseudocode,algorithm,wasysym,epstopdf,verbatim,array,textcomp,balance}
\usepackage[table]{xcolor}
\DeclareMathOperator{\sign}{sign}
\DeclareMathOperator{\nullspace}{null}

\DeclareMathOperator{\diag}{dg}

\newtheorem{assumption}{Assumption}
\newtheorem{proposition}{Proposition}
\newtheorem{lemma}{Lemma}

\newtheorem{theorem}{Theorem}
\newtheorem{definition}{Definition}

\newcommand \bzero{\mathbf{0}}
\newcommand \bone{\mathbf{1}}
\newcommand \ba{\mathbf{a}}

\newcommand \bd{\mathbf{d}}
\newcommand \be{\mathbf{e}}
\newcommand \bef{\mathbf{f}}
\newcommand \bg{\mathbf{g}} % exception since \bf is taken
\newcommand \bh{\mathbf{h}}

\newcommand \bn{\mathbf{n}}

\newcommand \bx{\mathbf{x}}

\newcommand \bA{\mathbf{A}}

\newcommand \bI{\mathbf{I}}

\newcommand \bS{\mathbf{S}}

\newcommand \bxi{\boldsymbol{\xi}}

\newcommand \tdf{\tilde{f}}
\newcommand \tdh{\tilde{h}}

\newcommand \tdn{\tilde{n}}

\newcommand \mcC{\mathcal{C}}

\newcommand \mcG{\mathcal{G}}

\newcommand \mcL{\mathcal{L}}

\newcommand \mcN{\mathcal{N}}

\newcommand \mcP{\mathcal{P}}

\newcommand \mcT{\mathcal{T}}

\newcommand \bmcP{\bar{\mathcal{P}}}

\newcommand \tbf{\tilde{\mathbf{f}}}
\newcommand \tbh{\tilde{\mathbf{h}}}

\newcommand \tbn{\tilde{\mathbf{n}}}

\newcommand \tbA{\tilde{\mathbf{A}}}

\begin{document}
\title{On the Flow Problem in Water Distribution Networks: Uniqueness and Solvers}
%\title{Solving the Flow Problem in Water Distribution Networks}
%\title{Water Flow Problem: Uniqueness and Novel Solvers}

\author{Manish K. Singh~\IEEEmembership{Student Member,~IEEE} and Vassilis Kekatos,~\IEEEmembership{Senior Member,~IEEE}

\thanks{Manuscript received February 12, 2019; revised August 13, 2019, March 31, 2020, and July 31, 2020; accepted September 17, 2020. Date of publication DATE; date of current version DATE. Paper no. TCONES-20-0170.}

\thanks{M. K. Singh and V. Kekatos are with the Bradley Dept. of Electrical and Computer Engineering, Virginia Tech, Blacksburg, VA 24061, USA. Emails: \{manishks,kekatos\}@vt.edu. This work was supported by the U.S. National Science Foundation under Grant 1711587.}
\thanks{Color versions of one or more of the figures is this paper are available online at {http://ieeexplore.ieee.org}.}
\thanks{Digital Object Identifier XXXXXX}
%\vspace*{-1.5em}
}	

\markboth{IEEE TRANSACTIONS ON CONTROL OF NETWORK SYSTEMS (to appear)}{Singh and Kekatos: On the Flow Problem in Water Distribution Networks: Uniqueness and Solvers}
		
\maketitle
	
\begin{abstract}
Increasing concerns on the security and quality of water distribution systems (WDS), call for computational tools with performance guarantees. To this end, this work revisits the physical laws governing water flow and provides a hierarchy of solvers of complementary value. Given the water injection or pressure at each WDS node, finding the water flows within pipes and pumps along with the pressures at all WDS nodes constitutes the water flow (WF) problem. The latter entails solving a set of (non)-linear equations. We extend uniqueness claims on the solution to the WF equations in setups with multiple fixed-pressure nodes and detailed pump models. For networks without pumps, the WF solution is already known to be the minimizer of a convex function. The latter approach is extended to networks with pumps but not in cycles, through a stitching algorithm. For networks with non-overlapping cycles, a provably exact convex relaxation of the pressure drop equations yields a mixed-integer quadratically-constrained quadratic program (MI-QCQP) solver. A hybrid scheme combining the MI-QCQP with the stitching algorithm can handle WDS with overlapping cycles, but without pumps on them. Each solver is guaranteed to converge regardless of initialization, as numerically validated on a benchmark WDS. 
\end{abstract}
	
	
%%%%%%%%%%%%%%%%%%%%%%%%%%%%%%%%%%%%%%%%%%%%%%%%%%%%%%%%%%%%%%%%%%
\section{Introduction}\label{sec:intro}
\allowdisplaybreaks
Water distribution systems serve as a critical infrastructure across the world. The direct dependence of human lives on the availability of water has motivated research on the security, resiliency, and quality of water supply systems~\cite{Helena2017lost}, \cite{Rossman94net2}. The high cost of installation for different WDS components renders long-term network planning an important problem~\cite{Sherali2001}, \cite{bragalli2015mathprog}. Furthermore, the relatively expensive operation of a WDS is primarily attributed to the electricity cost for running pumps to properly circulate water~\cite{Vanzyl2004}. Thus, optimal pump scheduling for the daily operation of a WDS is a pertinent research problem~\cite{TaylorOWFCDC}, \cite{TaylorOWFcones}, \cite{singh2018optimal}, \cite{taha2019accgp}. An inevitable component of the aforementioned computational problems is satisfying the physical laws governing water flow. Mathematically, the \emph{water flow (WF) equations} consist of a set of linear equations ensuring mass conservation, along with a set of non-linear equations arising from energy and momentum conservation~\cite{Rossman2013unified}.

Solving the WF equations constitutes the \emph{water flow problem}. Specifically, given water demand at all nodes, the standard WF task aims at finding the water flows within all pipes and the pressures at all nodes complying with the WF equations. Despite nonlinear, these equations enjoy a unique solution for networks with a single fixed-pressure node and when the pressure added by pumps is approximated as constant~\cite{Humpola2013unified}. Modern renditions of the WF problem may incorporate pumps, valves, and pressure-based demands~\cite{jun2013pda}, \cite{epanet2000}. Either way, handling the non-linear equations stemming from the conservation of energy and momentum remain the core challenge~\cite{Rossman2013unified}. 

Existing WF solvers update iteratively a set of WF variables, which could be the pipe flows, loop flows, nodal pressures, or combinations thereof~\cite{Todini2006convergence}. These solvers can be broadly classified into those relying on successive linear approximations, and those relying on Newton-Raphson-type of updates~\cite{Rossman2013unified}. For example, the WF solver of \cite{Hafez18} constitutes a fixed-point iteration and belongs to the former class, while EPANET (perhaps the most widely used WF solver) to the latter~\cite{epanet2000}. In fact, most of the schemes within each class have been shown to be equivalent to each other upon a (non)-linear transformation of variables~\cite{Todini2006convergence},~\cite{Rossman2013unified}.  To cope with the problem of dimensionality, different preprocessing, partitioning, and reformulations have been reported in~\cite{elhay2014cotree}, \cite{vidal2015loop}, \cite{elhay2018partition}. 

The aforesaid solvers exhibit two major shortcomings. First, their convergence and rate of convergence depend critically on initialization. For example, it has been numerically demonstrated that EPANET fails to find a WF solution for some practical water networks~\cite{Zhang2017fixed}, \cite{estrada2009irrigation}. Nonetheless, proper initialization may be challenging when dealing with stochastic planning or risk analysis, where the WF task has to be solved repeatedly and under varying demands~\cite{vidal2015loop}. As a second shortcoming, the existing solvers do not naturally extend to optimal water flow (OWF) formulations. Therefore, most OWF efforts resort to linear approximations; non-linear local optimization; or slow zero-order algorithms building on an independent WF solver such as EPANET; see~\cite{Helena2017lost}. Such OWF approaches lack scalability and/or optimality guarantees. 

The contribution of this work is in two fronts. After a brief modeling of water networks (Section~\ref{sec:WDS}), this work first generalizes the uniqueness of the WF solution to WDS setups with multiple fixed-pressure nodes and more practical pump models with flow-dependent pressure gains. Second, it puts forth a suite of solvers that can provably recover the WF solution under different network setups; see also Fig.~\ref{fig:flowchart}:
\renewcommand{\labelenumi}{\emph{\roman{enumi})}}
\begin{enumerate}
\item In networks without pumps, it is already known that the WF solution can be recovered as the minimizer of a convex problem~\cite{Humpola2013unified}. Sections~\ref{subsec:cWF} devises SOCP- and dual decomposition-based solvers, while Section~\ref{subsec:hw} applies them to other pressure drop laws.
\item This energy function-based approach is extended to networks where pumps do not lie on cycles through the \emph{stitching algorithm} of Section~\ref{subsec:stitching}.
\item In networks with no overlapping cycles, the mixed-integer quadratically-constrained quadratic program (MI-QCQP) of Section~\ref{sec:MI-QCQP} can find the WF solution.
\item A hybrid procedure combining the stitching algorithm and the MI-QCQP solver handles networks having overlapping cycles, but without pumps on them (Section~\ref{sec:hybrid}).
\end{enumerate}
These novel solvers not only operate under different network configurations, but constitute a hierarchy: Solver \emph{ii)} builds on \emph{i)}; and solver \emph{iv)} builds on \emph{ii)} and \emph{iii)}. Numerical tests on benchmark networks evaluate the correctness and running times for \emph{i)} and \emph{iii)}; and validate them against the EPANET solver (Section~\ref{sec:tests}).
		
Regarding notation, lower- (upper-) case boldface letters denote column vectors (matrices). Calligraphic symbols are reserved for sets. The vectors of all zeros, all ones, and the $n$-th canonical vector are denoted respectively by $\bzero$, $\bone$, and $\be_n$, while their dimension will be clear from the context. The symbol $^{\top}$ stands for transposition.

%%%%%%%%%%%%%%%%%%%%%%%%%%%%%%%%%%%%%%%%%%%%%%%%%%%%%%%%%%%%%%
\section{Water Distribution System Modeling}\label{sec:WDS}
A WDS can be represented by a directed graph ${\mcG:=(\mcN,\mcP)}$. Its nodes are indexed by $n \in \mcN:=\{1,\dots,N\}$ and correspond to water reservoirs, tanks, and points of water demand. Let $d_n$ be the rate of water injected into the WDS from node $n$. For reservoirs apparently $d_n\geq 0$; for nodes with water consumers $d_n\leq0$; tanks may be filling or emptying; and $d_n=0$ for junction nodes. 

The edges in set $\mcP$ with cardinality $P:=|\mcP|$, are associated with pipes and pumps. The directed edge $p=(m,n)\in\mcP$ models the water pipe between nodes $m$ and $n$. Its water flow is denoted by $f_{mn}$ or $f_p$ depending on the context. If water flows from node $m$ to $n$, then $f_{mn}\geq0$; otherwise $f_{mn}<0$. 

Conservation of water flow dictates that for all $n\in\mcN$
\begin{equation*}
d_n=\sum_{k:(n,k)\in\mcP} f_{nk} - \sum_{k:(k,n)\in\mcP} f_{kn}.
\end{equation*}
The connectivity of the WDS is captured by the edge-node incidence matrix $\bA\in\mathbb{R}^{P\times N}$ with entries
	\begin{equation}\label{eq:Amatrix}
	A_{p,k}=
	\begin{cases}
	+1, &k=m\\
	-1, &k=n\\
	0, &\text{otherwise}
	\end{cases}~\forall~p=(m,n)\in\mcP.
	\end{equation}
Given $\bA$ and upon stacking flows and injections respectively in $\bef\in\mathbb{R}^{P}$ and $\bd\in\mathbb{R}^N$, the conservation of water flow across the WDS can be compactly expressed as
	\begin{equation}\label{eq:wmc}
	\bA^\top\bef=\bd.
	\end{equation}

	The operation of WDS is also governed by pressures. Water pressure is surrogated by \emph{pressure head}, defined as the equivalent height of a water column in meters, which exerts the surrogated pressure at its bottom. The pressure heads at all WDS nodes are measured with respect to a common geographical elevation level. The pressure head (henceforth \emph{pressure}) at node $n$ is denoted by $h_n$. Moreover, the pressure $h_r$ at a reservoir node $r\in\mcN$ typically serves as the pressure of reference. 
	
	With water flowing in a pipe, pressure drops along the direction of flow due to friction. The pressure drop across pipe $(m,n)\in \mcP$ is described by the Darcy-Weisbach law~\cite{Verleye_tankpumpmodel}
	\begin{equation}\label{eq:headloss}
	h_m-h_n=c_{mn}\sign(f_{mn})f_{mn}^2
	\end{equation}
where the constant $c_{mn}$ depends on pipe dimensions~\cite{singh2018optimal}, and the sign function is defined as 
\[\sign(x):=\left\{\begin{array}{rl}
+1,&~x>0\\
0,&~x=0\\
-1,&~x<0 
\end{array}\right..\]
The pressures at all nodes are collected in vector $\bh\in\mathbb{R}^N$. 
	
	%Adding this common minimum value of manometric pressure to the specific but known geographical elevation of each node $m\in\mcM$ gives a lower limit on its pressure as
	%\begin{equation}\label{eq:hmin}
	%h_m\geq \underline{h}_m.
	%\end{equation} 
	
To maintain pressures at desirable levels, water utilities use pumps on specific pipes. Let $\mcP_a\subset\mcP$ be the subset of pipes hosting a pump. The pipes in $\mcP_a$ can be considered lossless; this is without loss of generality since a pump can be modeled by an ideal pump followed by a short pipe. The remaining edges form the subset $\bmcP_a:=\mcP\setminus\mcP_a$, and correspond to lossy pipes governed by \eqref{eq:headloss}. When pump $p=(m,n)\in \mcP_a$ is running, it adds pressure $g_{mn}\geq 0$ so that
\begin{equation*}%\label{eq:headgain}
h_n-h_m=g_{mn}.
\end{equation*}
As detailed below, the pressure added by a pump decreases with the water flow through the pump \cite{CohenQHmodel}, \cite{Ulanicki_pumpmodel}. Moreover, for variable speed pumps, the pressure added increases with the pump speed. The exact relation is provided by manufacturers in the form of pump operation curves, and are oftentimes approximated with quadratic curve fits~\cite{CohenQHmodel}, \cite{Ulanicki_pumpmodel}, \cite{TaylorOWFcones}. In detail, the pressure added by pump $p=(m,n)$ is modeled as
\begin{equation}\label{eq:pumphead0}
g_{mn}(f_{mn},\omega_{mn})=\lambda_{mn} f_{mn}^2+\mu_{mn}\omega_{mn} f_{mn} + \nu_{mn} \omega_{mn}^2
\end{equation}
where $\omega_{mn}$ is the pump speed; and $\lambda_{mn}<0$, $\mu_{mn}\geq0$, and $\nu_{mn}\geq0$ are known pump parameters~\cite{TaylorOWFcones}. When pump $(m,n)$ is running, its flow has to be maintained within the range $0\leq\underline{f}_{mn}\leq f_{mn}\leq\bar{f}_{mn}$ due to engineering limitations~\cite{Verleye_tankpumpmodel}. Whereas the speed for fixed-speed pumps is constant; for variable-speed pumps, it is controlled by the operator. Either way, for a WF task the speed of every pump is given; otherwise the problem would become under-determined. Given its speed $\omega_{mn}^0$, the pressure added by pump $(m,n)$ is
\begin{equation}\label{eq:pumphead}
g_{mn}(f_{mn};\omega^0_{mn})=h_n-h_m=\lambda_{mn} f_{mn}^2+\bar{\mu}_{mn} f_{mn} + \bar{\nu}_{mn}
\end{equation}
where $\bar{\mu}_{mn}:=\mu_{mn}\omega_{mn}^0$ and $\bar{\nu}_{mn}:=\nu_{mn}(\omega_{mn}^0)^2$. Because the pump parameters satisfy $2\lambda_{mn} \underline{f}_{mn} + \mu_{mn}\omega_{mn}<0$ over the operating range of pump speeds $\omega_{mn}$, the pressure gain due to any pump is a \emph{strictly decreasing function} of the water flow in the range $[\underline{f}_{mn},\bar{f}_{mn}]$. This observation is instrumental in establishing the uniqueness of the WF solution in Section~\ref{sec:uniqueness}.

When pump $p=(m,n)\in\mcP_a$ is not running, water can flow freely in either directions through a bypass valve~\cite{CohenQHmodel}. Because bypass valve sections are typically short, one can ignore the pressure drop along them to get $h_m=h_n$ and $g_{mn}=0$. In this case, the WDS graph can be reduced by removing pipe $p$ and node $n$, and connecting to node $m$ the edges previously incident to $n$. Alternatively, the valve can be modeled as a short lossy pipe, whose pressure drop is governed by \eqref{eq:headloss}.

Valves constitute a vital component for water flow control. They can be modeled by an on/off switch; a linear pressure-reducing model; a flow-dependent non-linear model; or a flow control model~\cite{TaylorOWFcones},~\cite{epanet2000}. Under the typical operational setup, the valve at the reference node regulates pressure, whereas the valves at the remaining reservoirs and tanks regulate flows.

Summarizing this section, a WDS operating point is described by the triplet $(\bd,\bef,\bh)$ satisfying the WF equations of \eqref{eq:wmc},~\eqref{eq:headloss}, and \eqref{eq:pumphead}. This work deals with the uniqueness of a WF solution and efficient solvers for finding this solution.

%%%%%%%%%%%%%%%%%%%%%%%%%%%%%%%%%%%%%%%%%%%%%%%%%%%%%%%%%%%%%%
\section{Uniqueness of Water Flow Solution}\label{sec:uniqueness}
Different from the OWF problem where tanks and pumps are scheduled over a time horizon, the WF task aims at solving the WF equations given the water injections $\bd$. As such, it constitutes a key component of WDS operation and planning. The WF problem is formally stated next. 

\begin{definition}\label{def:WF}
Given: i) water injections $\bd$; ii) the statuses and speeds $\{\omega_{mn}\}$ for all pumps $(m,n)\in\mcP_a$; and iii) the pressure $h_r$ at the reference node $r\in\mcN$; the WF task aims at finding the flows $\bef$ and pressures $\bh$ satisfying \eqref{eq:wmc},~\eqref{eq:headloss}, \eqref{eq:pumphead}.
\end{definition}

The WF task involves $N+P-1$ equations over $N+P-1$ unknowns. The water balance in \eqref{eq:wmc} yields $N-1$ linearly independent equations. In addition, the pressure drops across lossy pipes [cf.~\eqref{eq:headloss}], and the pressure gains due to pumps [cf.~\eqref{eq:pumphead}] provide $P$ non-linear equations.

According to Definition~\ref{def:WF}, both $\bef$ and $\bh$ are unknown. Nevertheless, to find a triplet $(\bd,\bef,\bh)$ that satisfies \eqref{eq:wmc},~\eqref{eq:headloss}, and \eqref{eq:pumphead}, it suffices to find either $\bef$ or $\bh$. To see this, note that if $\bef$ is known, then $\bh$ can be calculated from \eqref{eq:headloss}--\eqref{eq:pumphead} and $h_r$. On the other hand, if $\bh$ is known, the flows $\bef$ can be found thanks to the monotonicity of \eqref{eq:headloss} and the monotonicity of \eqref{eq:pumphead} in $[\underline{f}_{mn},\bar{f}_{mn}]$. In a nutshell, solving the WF task amounts to finding either $\bef$ or $\bh$. Because this simple observation is used throughout our analysis, it is summarized as a lemma.

\begin{lemma}\label{le:fh}
Given $\bd$, a triplet $(\bd,\bef,\bh)$ satisfying \eqref{eq:wmc},~\eqref{eq:headloss}, and \eqref{eq:pumphead} is uniquely characterized by $\bef$ or $\bh$.
\end{lemma}

The WF task can be posed as the feasibility problem
\begin{align*}\label{eq:W1}
	\mathrm{find}~&~\{\bef,\bh\}\tag{W1}\\
	\mathrm{s.to}~&~ \eqref{eq:wmc},\eqref{eq:headloss},\eqref{eq:pumphead}.
\end{align*}
Since \eqref{eq:headloss}--\eqref{eq:pumphead} are quadratic equalities, problem~\eqref{eq:W1} is non-convex. For a tree WDS graph $\mcG$ (for which $P=N-1$), the edge-node incidence matrix $\bA$ has $N-1$ linearly independent rows~\cite{GodsilRoyle}. Then the flows $\bef$ can be found uniquely from \eqref{eq:wmc}, and the WF task is readily solved according to Lemma~\ref{le:fh}. However, handling the WF task in a loopy $\mcG$ containing pumps remains non-trivial. Naturally, analyzing the uniqueness of a WF solution is a critical task. Reference~\cite{Humpola2013unified} establishes the uniqueness of a WF solution under the assumption that the pressure gain added by a pump is constant. The next claim extends this uniqueness result regardless of the structure of $\mcG$ and to the more detailed pump model of ~\eqref{eq:pumphead}.
	
\begin{theorem}\label{th:unique}
If the WF equations are feasible for some $\bd$, they feature a unique solution.
\end{theorem}
	
\begin{IEEEproof}
Proving by contradiction, assume $(\bh,\bef)$ and $(\tbh,\tbf)$ are two distinct solutions of \eqref{eq:W1}. Since the two flow vectors $\bef$ and $\tbf$ satisfy~\eqref{eq:wmc}, the vector $\bn:=\tbf-\bef$ lies in the nullspace of $\bA^\top$, that is $\bA^\top\bn=\bzero$. Then it follows that
		\begin{equation}\label{eq:h12}
		\bn^\top\bA(\tbh-\bh)=0.
		\end{equation}
		Let us decompose $\bn$ into its positive and negative entries as $\bn=\bn_+-\bn_-$, where $\bn_+\geq \bzero$ and $\bn_-\geq \bzero$. Plugging this decomposition into \eqref{eq:h12} yields
		\begin{equation}\label{eq:thwfuni1}
		\bn_+^\top\left(\bA\tbh-\bA\bh\right)=\bn_-^\top\left(\bA\tbh-\bA\bh\right).
		\end{equation}
		%The LHS in \eqref{eq:thwfuni1} is a weighted sum of changes in pressure drops while transitioning from $\bef_1$ to $\bef_2$, across edges corresponding to positive entries of $\bn$. The evaluation in the RHS of \eqref{eq:thwfuni1} similarly relates to the edges corresponding to negative entries of $\bn$.
		
Recall from \eqref{eq:headloss} that the pressure drop across a lossy pipe is monotonically increasing in flow. Similarly, the pressure added by a pump described by \eqref{eq:pumphead} is monotonically decreasing in flow. In other words, the pressure drop along a pump is monotonically increasing in flow. Hence, for each edge $p\in\mcP$:
		\begin{equation}\label{eq:thwfuni2}
		\ba_p^\top \tbh>\ba_p^\top \bh \quad \textrm{if and only if} \quad \tilde{f}_p > f_p 
		\end{equation}
where $\ba_p^\top$ is the $p$-th row of matrix $\bA$.
		
Consider the $p$-th entry of vector $\bn$. If $\tilde{f}_p>f_p$, then $n_{+,p}>0$ and $n_{-,p}=0$ by definition of $\bn_+$ and $\bn_{-}$. Then edge $p$ contributes to the left-hand side (LHS) of \eqref{eq:thwfuni1}. The monotonicity in \eqref{eq:thwfuni2} entails that $\ba_p^\top (\tbh-\bh)>0$ and so $n_{+,p}\cdot\ba_p^\top (\tbh-\bh)>0$. The latter holds for all edges contributing to the LHS of \eqref{eq:thwfuni2}, so that 
		\[\bn_+^\top\left(\bA\tbh-\bA\bh\right)>0.\]
		
		On the other hand, if $\tilde{f}_p<f_p$, then $n_{+,p}=0$ and $n_{-,p}>0$. Then, edge $p$ contributes to the right-hand side (RHS) of \eqref{eq:thwfuni2}. The monotonicity in \eqref{eq:thwfuni2} entails $n_{-,p}\cdot \ba_p^\top (\tbh-\bh)<0$. Applying the claim over all edges participating in the RHS of \eqref{eq:thwfuni2} provides
		\[\bn_-^\top\left(\bA\tbh-\bA\bh\right)<0.\] 
		The signs of the LHS and RHS contradict the equality in \eqref{eq:thwfuni1}, thus proving the claim.
	\end{IEEEproof}

Regarding the pressure drop law of \eqref{eq:headloss}, the Hazen-Williams equation is sometimes used wherein the flow $f_{mn}$ is raised to the exponent of $1.852$. This exponent is different from the exponent of $2$ in the Darcy-Weisbach equation; see for example~\cite{Hafez18}. While the Darcy-Weisbach equation is a theoretical formula, the Hazen-Williams equation is based on curve fitting of experimental data~\cite{bragalli2015mathprog}. Either way, the uniqueness argument of Theorem~\ref{th:unique} holds for any positive exponent on the water flow $f_{mn}$ involved in the pressure drop equation of \eqref{eq:headloss}.

The WF problem posed in Definition~\ref{def:WF} follows the most traditional setup pursued in WF literature~\cite{Humpola2013unified},~\cite{Collins78}, \cite{Maugis77}. However, with multiple sources feeding a WDS an alternate form of WF becomes pertinent. Specifically, multiple supply nodes may operate as fixed-pressure nodes while the remaining nodes act as fixed-injection nodes~\cite{taha2020WFGP}. The ensuing result extends the uniqueness of a WF solution to this more general setting. 

\begin{proposition}\label{le:unique2}
Given: i) pressures $h_i$ for $i$ in a subset of nodes $\mcN_r\subset\mcN$; ii) injections $d_i$ for all $i\in\mcN\setminus\mcN_r$; and iii) the statuses and speeds for all pumps; the set of equations \eqref{eq:wmc},~\eqref{eq:headloss}, and \eqref{eq:pumphead}  has at most one solution.
\end{proposition}

\begin{IEEEproof}[Proof outline] Consider two pairs $(\bd,\bef)$ and $(\bd',\bef')$ satisfying \eqref{eq:wmc}. If $\bd=\bd'$, the uniqueness of the WF solution follows from Theorem~\ref{th:unique}. Otherwise, \cite[Lemma~3]{singh2018TCNSGF} dictates that there exists a path $\mcP_{mn}$ between fixed-pressure nodes $m$ and $n$, along which the flows in $\bef$ differ from those in $\bef'$ in a consistent direction, i.e., either $f_p>f_p'$ or $f_p<f_p'$ for all edges $p\in\mcP_{mn}$. However, in that case, the monotonicity argument used in \eqref{eq:thwfuni2} would entail for the respective pressures $h_m-h_n\neq h_m'-h_n'$, which contradicts  the hypothesis of $m$ and $n$ being fixed-pressure nodes.
\end{IEEEproof}

\begin{figure}[t]
	\centering
	\includegraphics[scale=0.43]{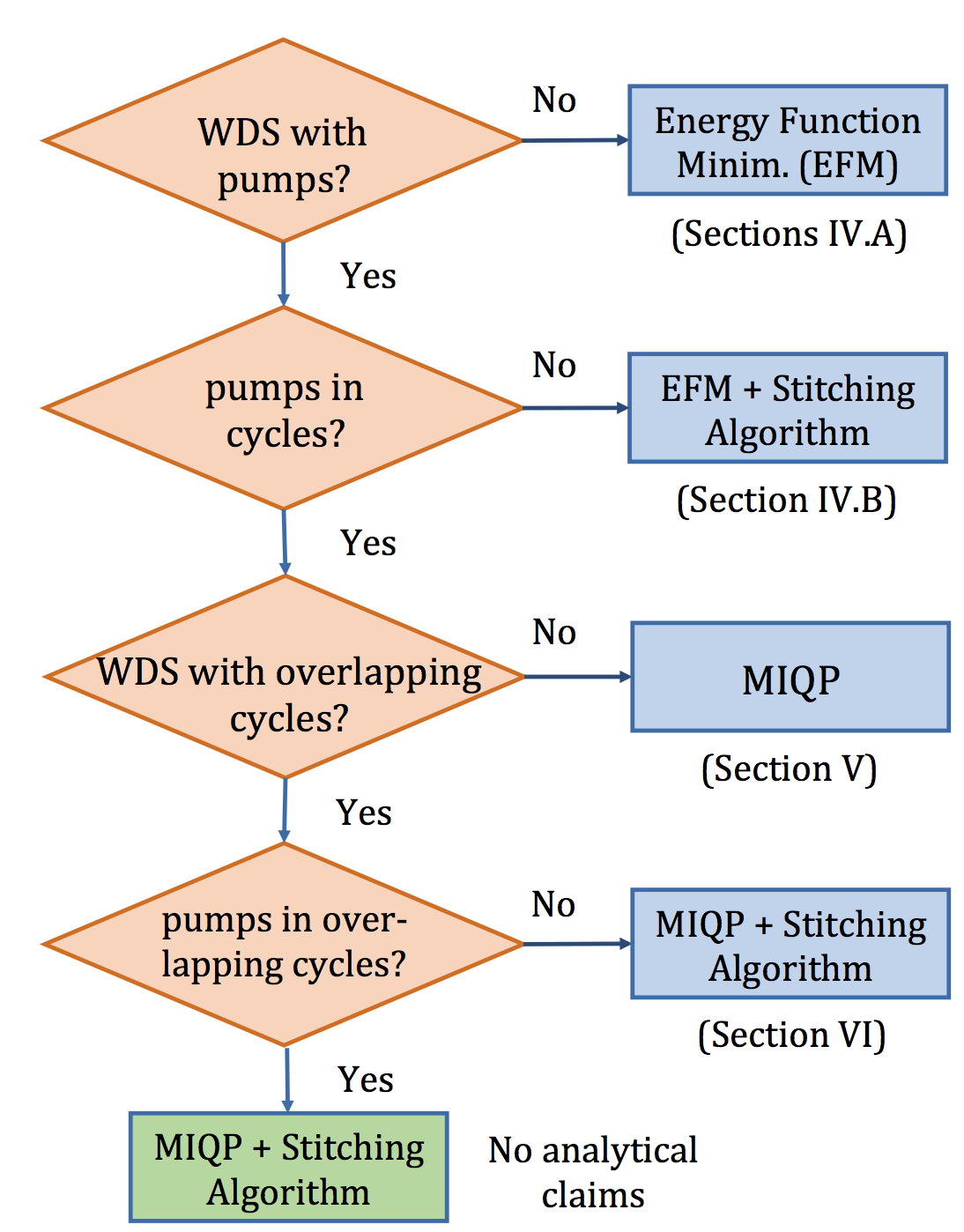}
	\caption{Given a WDS, this flowchart suggests the most suitable WF solver based on its analytical guarantees and computational complexity.}
	\label{fig:flowchart}
\end{figure}

Having established the uniqueness of the WF solution, the ensuing sections develop a suite of WF solvers of complementary value: Each solver finds provably the WF solution under different network setups. The devised solvers exhibit different complexity, yet none of them requires a proper initialization. Figure~\ref{fig:flowchart} summarizes the particular WDS setups each solver can handle, and serves as a roadmap for the following sections. The solvers developed in this work are for the WF task as posed in Definition~\ref{def:WF}, whereas WF solvers for the setting considered in Proposition~\ref{le:unique2} are beyond the scope of this work. 

%\cmb{To this end, three approaches are next presented which are useful for different network conditions. The first approach relies on solving a convex constrained minimization problem, and has been adapted from the literature on solving gas flow problems. The second approach is an alternative to the first, where an unconstrained convex minimization problem has been formulated which yields the pressure vector $\bh$. Finally, a mixed-integer penalized relaxation is presented that solves the WF problem for a much larger class of water networks.}

%%%%%%%%%%%%%%%%%%%%%%%%%%%%%%%%%%%%%%%%%%%%%%%%%%%%%%%%%%%%%%	
\section{Energy Function-Based Water Flow Solvers}\label{sec:energy}
For a WDS without pumps, the WF task simplifies to solving \eqref{eq:wmc}--\eqref{eq:headloss}. These equations are structurally similar to the equations governing the flows and pressures in a natural gas network under steady-state conditions~\cite{singh2018TCNSGF},~\cite{Singh18ACC}. Reference~\cite{Collins78} finds the solution to the nonlinear WF equations in networks without pumps as the minimizers of a convex \emph{energy function}; see also~\cite{Maugis77} and \cite{wolf2000energy} for counterparts in gas networks. Upon briefly reviewing these reformulations, this section expands their scope to WDS with no pumps in cycles and other pressure drop laws; see roadmap of Figure~\ref{fig:flowchart}. %The solver of Section~\ref{subsec:stitching} in particular serves as the foundation for the hybrid solvers of Section~\ref{sec:hybrid}.

\subsection{Finding the WF Solution in WDS without Pumps}\label{subsec:cWF}
In the absence of pumps, the WF problem \eqref{eq:W1} has been posed as a constrained minimization over $\bef$~\cite{Collins78}.

\begin{lemma}[\cite{Collins78}]\label{le:efm}
In a WDS without pumps ($\mcP=\bmcP_a$), the vector of water flows $\bef$ satisfying \eqref{eq:wmc}--\eqref{eq:headloss} can be found as the unique minimizer of the convex minimization problem
\begin{subequations}\label{eq:cWF}
\begin{align}
\min_{\bef}~&~\frac{1}{3}\sum_{p\in\mcP} c_p|f_p|^3 \label{eq:constrainedWF:cost}\\
\mathrm{s.to}~&~\mathbf{A}^{\top} \bef = \bd. \label{eq:constrainedWF:con}
\end{align}
\end{subequations}
Moreover, if $\bxi^*$ is the vector of optimal Lagrange multipliers corresponding to \eqref{eq:constrainedWF:con}, then the nodal pressures are provided by $\bh=(\bI - \bone\be_r^\top)\bxi^* +h_r\bone$. 
\end{lemma}

%\begin{IEEEproof}
%%The objective in \eqref{eq:cWF} can be written as the $\ell_3$-norm of a scaled flow vector raised to the third power. Its strict convexity follows from composition rules~\cite[Sec.~3.2.4]{BoVa04}. Then, the uniqueness of the minimizer of \eqref{eq:cWF} stems from strong convexity. Moreover, the objective is differentiable: the $p$-th entry of its gradient vector $\bg_c(\bef)$ is $c_p\sign(f_p) |f_p|^2$.
%%Moreover, its Hessian matrix is diagonal, and its $(p,p)$-th entry equals $2c_p |f_p|$. The positive semi-definiteness of this Hessian matrix entails that the objective function is convex.
%The Lagrangian function associated with \eqref{eq:cWF} is
%\begin{equation}\label{eq:Lagrangian}
%L(\bef;\bxi)=\frac{1}{3}\sum_{p\in\mcP} c_p |f_p|^3+\bxi^\top(\bd-\bA^\top \bef).
%\end{equation}
%The optimal primal and dual variables $(\bef^*,\bxi^*)$ satisfy the Lagrangian optimality condition, that is $\bg_c(\bef^*)=\bA\bxi^*$, which coincides with \eqref{eq:headloss}. By the definition of $\bA$ in \eqref{eq:Amatrix}, it holds that $\bA\bxi^*=\bA(\bxi^*+\alpha\bone)$ for any $\alpha\in\mathbb{R}$. Hence, the pressures can be recovered either by Lemma~\ref{le:fh}, or by fixing $\alpha=h_r-\be_r^\top\bxi^*$ to equate the $r$-th entry of $\bxi^*+\alpha\bone$ to the reference pressure $h_r$.
%\end{IEEEproof}

The claim follows from the optimality conditions of \eqref{eq:cWF}. The strict convexity of \eqref{eq:constrainedWF:cost} shows that the WF solution is unique in networks \emph{without pumps}. This claim has also been established using a contraction argument in \cite{Hafez18}. Theorem~\ref{th:unique} and Proposition~\ref{le:unique2} generalize this uniqueness claim to networks having pumps and multiple fixed-pressure nodes.

 We next present two ways for handling \eqref{eq:cWF}: an second-order cone program (SOCP) reformulation, and a dual decomposition approach. It is known that  $\ell_p$-norms can be handled using a hierarchy of second-order cone (SOC) constraints; see~\cite[Eq.~(11)]{alizadeh2003second}. In the context of WDS, a more compact hierarchy of SOCs was originally put forth in \cite{TaylorOWFCDC} to handle the cubic powers of \eqref{eq:constrainedWF:cost} involved in an OWF problem. Adopting the latter reformulation, the minimizer of \eqref{eq:cWF} can be found by solving the SOCP
\begin{subequations}\label{eq:socp}
\begin{align}
\min_{\{f_p,w_p,y_p,t_p\}}~&~\frac{1}{3}\sum_{p\in\mcP} c_p t_p \label{eq:socp:cost}\\
\mathrm{s.to}~&~\bA^{\top} \bef = \bd. \label{eq:socp:con1}\\
&~-w_p\leq f_p\leq w_p,&  \forall p\label{eq:socp:con2}\\
&~w_p^2\leq y_p,~y_p^2\leq w_p t_p, &\forall p.\label{eq:socp:con3}
\end{align}
\end{subequations}
The fact that the $f_p$'s minimizing \eqref{eq:cWF} and \eqref{eq:socp} coincide can be established by showing that the solution of \eqref{eq:socp} satisfies $w_p=|f_p|$; $y_p=w_p^2=f_p^2$; and $t_p=w_p^3=|f_p|^3$ for all $p\in\mcP$; see also \cite{TaylorOWFCDC}--\cite{TaylorOWFcones} for details. Each constraint in \eqref{eq:socp:con3} is a rotated second-order cone~\cite{alizadeh2003second}.

%In the OWF approaches of \cite{TaylorOWFCDC}--\cite{TaylorOWFcones}, flow directions are presumed known (say $f_{mn}\geq 0$) and the cost of \eqref{eq:constrainedWF:cost} is minimized subject to \eqref{eq:headloss} and other OWF constraints. Equation~\eqref{eq:headloss} reading now $h_m-h_n=c_{mn}f_{mn}^2$ is relaxed to $h_m-h_n\geq c_{mn}f_{mn}^2$ to yield a mixed-integer conic program, yet the constraint turns out to be satisfied with equality under certain conditions. In contrast to \cite{TaylorOWFCDC}--\cite{TaylorOWFcones}, the constraints \eqref{eq:socp:con2}--\eqref{eq:socp:con3} capture the $|f_p|^3$ terms of \eqref{eq:constrainedWF:cost}, and no relaxation of the Darcy-Weisbach law. Further, problems \eqref{eq:cWF} and \eqref{eq:socp} are equivalent; not a relaxation of one another; do not even involve pressures as variables; and $\bd$ is given. \cmr{This para may be trimmed/removed?}

Alternatively, problem \eqref{eq:cWF} can be tackled via dual decomposition: During its $k$-th iteration, the primal variable $\bef$ is updated by minimizing the Lagrangian function $L(\bef;\bxi^k)$ evaluated at the latest estimate of dual variables $\bxi^k$. The latter problem decouples across pipes and enjoys a closed-form solution as
\begin{equation*}
f_p^{k+1}:=\arg\min_{f_p}\frac{c_p|f_p|^3}{3}-f_p\ba_p^\top \bxi^k=\sign(\ba_p^\top \bxi^k)\sqrt{\frac{|\ba_p^\top \bxi^k|}{c_p}}.
\end{equation*}
For a step-size $\mu>0$, the dual variables are updated as
\[\bxi^{k+1}:=\bxi^k+\mu\left(\bd-\bA^\top\bef^{k+1}\right).\]

%%%%%%%%%%%%%%%%%%%%%%%%%%%%%%%%%%%%%%%%%%%%%%%%%%%%%%%%%%%
%\subsection{Unconstrained Formulation}\label{subsec:uWF}

Again for WDS without pumps, the WF task can alternatively be posed as an \emph{unconstrained} minimization over $\bh$~\cite{Collins78}
%\begin{equation}\label{eq:uWF}
%\min_{\bh}~\frac{1}{3}\sum_{n=1}^N\sum_{m:m\sim n}\frac{|h_m-h_n|^\frac{3}{2}}{\sqrt{c_{mn}}} - \bd^\top\bh
%\end{equation}
%where the notation $m\sim n$ means that nodes $m$ and $n$ are connected, that is $(m,n)\in\mcP$ or $(n,m)\in\mcP$. Note that each edge $(m,n)\in\mcP$ contributes twice in the summation of \eqref{eq:uWF}.
\begin{equation}\label{eq:uWF}
\min_{\bh}~\frac{2}{3}\sum_{(m,n)\in\mcP}\frac{|h_m-h_n|^\frac{3}{2}}{\sqrt{c_{mn}}} - \bd^\top\bh.
\end{equation}

The objective of \eqref{eq:uWF} is convex (by composition rules) and differentiable. The $n$-th entry of its gradient $\bg_u(\bh)$ is
\[[\bg_u(\bh)]_n=\sum_{p=(m,n)\in\mcP}\sign(\ba_p^\top \bh)\sqrt{\frac{|\ba_p^\top \bh|}{c_p}}-d_n.\]
Setting this gradient equal to zero yields the WF equations after eliminating $\bef$ from \eqref{eq:wmc}--\eqref{eq:headloss}. The ambiguity in pressures can be waived by shifting the minimizer of \eqref{eq:uWF} to match the reference pressure $h_r$. Once the pressures $\bh$ are found, the flow vector $\bef$ can be retrieved by Lemma~\ref{le:fh}. Problem \eqref{eq:uWF} is amenable to any first-order method for unconstrained optimization, such as the gradient descent iterations $\bh^{k+1}:=\bh^k-\mu\bg_u(\bh^k)$ for a step-size $\mu>0$, or accelerated variants.

\subsection{Extension to WDS with no Pumps in Cycles}\label{subsec:stitching}
Since problems \eqref{eq:cWF} and \eqref{eq:uWF} cannot handle water networks with pumps, this section extends their applicability to networks with pumps, but not on cycles. This is accomplished by adopting the reduction technique of \cite{mercado2002reduce} to build what we term \emph{stitching algorithm}:
	\renewcommand{\labelenumi}{\emph{S\arabic{enumi})}}
	\begin{enumerate}
		\item Remove the edges of $\mcG$ corresponding to pumps $\mcP_a$. The obtained graph contains $|\mcP_a|+1$ disconnected components $\mcG_c$ for $c=1,\ldots,|\mcP_a|+1$. 
		\item Replace each component $\mcG_c$ by a supernode, and connect the supernodes using the edges in $\mcP_a$ to create the supergraph $\mcG'$. Graph $\mcG'$ features a \emph{tree structure}.
		\item Find the total water injection per supernode $\mcG_c$. Since $\mcG'$ is a tree, the water flows on $\mcP_a$ can be found readily. 		
		\item If the flow along pump $(m,n)\in\mcP_a$ is $f_{mn}$, modify the injections at nodes $m$ and $n$ as $\hat{d}_m=d_m-f_{mn}$ and $\hat{d}_n=d_n+f_{mn}$.
%		\item Removing the pumps from the network results in a set of disconnected subnetworks. Nevertheless, water injections are known for every node.
		\item Solve \eqref{eq:cWF} per connected component $\mcG_c$ to find the water flows at all lossy pipes. 
		\item Having acquired vector $\bef$, the pressure vector $\bh$ can be found from Lemma~\ref{le:fh}.
\end{enumerate}

In \emph{S5)}, rather than solving the constrained minimization of \eqref{eq:cWF}, one could use its unconstrained counterpart of \eqref{eq:uWF}. Steps \emph{S1)}--\emph{S4)} are still needed to find a meaningful water injection vector per component $\mcG_c$. Once a pressure vector has been found per component, the pressures must be revised as follows: The pressures within the component containing the reference node, say component $\mcG_1$, are kept unaltered. Consider a pump running from node  $m\in\mcG_1$ to node $n\in\mcG_2$. Using the flow $f_{mn}$ computed in step \emph{S3)}, the pressure gain $g_{mn}=g_{mn}(f_{mn};\omega^0_{mn})$ can be found from \eqref{eq:pumphead}. Having found $h_m$ by solving \eqref{eq:uWF} in $\mcG_1$, the pressure at node $n$ can be computed as $h_n=h_m+g_{mn}$. If the pressure at node $n$ recovered by solving \eqref{eq:uWF} in $\mcG_2$ is $\tdh_n$, then all pressures within $\mcG_2$ should be shifted by $h_n-\tdh_n$. The process is repeated for all components to recover the entire pressure vector $\bh$. 

\subsection{Handling the Hazen-Williams Pressure Drop Law}\label{subsec:hw}
The aforementioned WF solvers can be modified to handle the Hazen-Williams in lieu of the Darcy-Weisbach pressure drop equations. Recall that according to the Hazen-Williams pressure drop law, the exponent of $|f_{mn}|$ in \eqref{eq:headloss} changes from $2$ to $\rho=1.852$. To handle this, problem~\eqref{eq:cWF} is generalized to
	\begin{subequations}\label{eq:cWF-HW}
		\begin{align}
		\min_{\bef}~&~\frac{1}{\rho+1}\sum_{p\in\mcP} c_p|f_p|^{\rho+1} \\
		\mathrm{s.to}~&~\mathbf{A}^{\top} \bef = \bd
		\end{align}
	\end{subequations}
where the exponent $\rho$ depends on the pressure drop law and could vary per pipe. Problem \eqref{eq:cWF-HW} remains convex for $\rho>0$. 

Even though the SOCP model of \eqref{eq:socp} cannot be used to solve \eqref{eq:cWF-HW}, a different SOCP reformulation can be derived as long as $\rho$ is rational; see~\cite[pp.~12--13]{alizadeh2003second} for details. Moreover, the dual decomposition updates of Section~\ref{subsec:cWF} can be adapted to solve \eqref{eq:cWF-HW}. Similarly, the unconstrained minimization of~\eqref{eq:uWF} can be modified to
\begin{equation}\label{eq:uWF-HW}
\min_{\bh}~\frac{\rho}{\rho+1}\sum_{(m,n)\in\mcP}\frac{|h_m-h_n|^{\frac{\rho+1}{\rho}}}{c_{mn}^{1/\rho}} - \bd^\top\bh.
\end{equation}

Heed that the energy function-based solvers for \eqref{eq:W1} can not handle the WF setups with multiple fixed-pressure nodes described under Proposition~\ref{le:unique2}. This is because energy function-based solvers do not take reference pressure as an input to the problem. Rather, the obtained pressures are shifted \emph{a posteriori} to agree with a given fixed pressure. For a single fixed-pressure node, the applicability of our WF solvers can be broadened to WDS with pumps in cycles. To do so, we next pursue an MI-QCQP solver. Unlike the energy function-based WF solvers of this section, the MI-QCQP solver of Section~\ref{sec:MI-QCQP} applies only to the Darcy-Weisbach equation; its possible extension to the Hazen-Williams pressure drop law goes beyond the scope of this work.

%%%%%%%%%%%%%%%%%%%%%%%%%%%%%%%%%%%%%%%%%%%%%%%%%%%%%%%%%%%%%%	
\section{MI-QCQP Water Flow Solver}\label{sec:MI-QCQP}
This section presents an MI-QCQP relaxation of \eqref{eq:W1} along with sufficient conditions for its exactness.
	
\subsection{Problem Reformulation}\label{subsec:PR}
The non-convexity of the WF problem \eqref{eq:W1} is due to the quadratic equality constraints of \eqref{eq:headloss} and \eqref{eq:pumphead}. Adopting the convex relaxation approach of \cite{TaylorOWFcones}--\cite{TaylorOWFCDC}, the pressure drop along pipe $(m,n)\in\bmcP_a$ is relaxed from \eqref{eq:headloss} to
	\begin{itemize}
		\item $h_m-h_n\geq c_{mn}f_{mn}^2$ for $f_{mn}\geq 0$; or
		\item $h_n-h_m\geq c_{mn}f_{mn}^2$ for $f_{mn}\leq 0$.
	\end{itemize}
To handle the two cases, let us introduce a binary variable $x_{mn}$ capturing the flow direction on pipe $(m,n)$. By using the so-termed big-$M$ trick, the two cases can be modeled as
	\begin{subequations}\label{eq:relaxed}
		\begin{align}
		-&M(1-x_{mn})\leq f_{mn}\leq Mx_{mn}\label{seq:hlra}\\
		-&M(1-x_{mn})\leq h_m-h_n-c_{mn}f_{mn}^2\label{seq:hlrb}\\
		&h_m-h_n + c_{mn}f_{mn}^2\leq Mx_{mn}\label{seq:hlrc}\\
		&x_{mn}\in\{0,1\}.\label{seq:hlrd}
		\end{align}
	\end{subequations}
The same MI-QCQP model has been also advocated for handling the pressure drop equations along the pipelines in natural gas networks; see e.g., \cite{backhaus2016convex},~\cite{bent2016hicss},~\cite{singh2018TCNSGF}. Constraint \eqref{seq:hlra} implies that $x_{mn}=\sign(f_{mn})$. If $x_{mn}=1$, the convex quadratic constraint in \eqref{seq:hlrb} is activated, whereas constraint \eqref{seq:hlrc} becomes trivial. The claim reverses for $x_{mn}=0$. Depending on the value $x_{mn}$, if either \eqref{seq:hlrb} or \eqref{seq:hlrc} is satisfied with equality, the relaxation is deemed as \emph{exact}. If that happens for all lossy pipes, the feasible set of \eqref{eq:relaxed} has captured the original non-convex constraints in \eqref{eq:headloss}. 

Similarly, the pressure added by pump $(m,n)\in\mcP_a$ is relaxed from \eqref{eq:pumphead} to
	\begin{equation}\label{eq:pumprel}
	h_n-h_m\leq\lambda_{mn}f_{mn}^2+\hat{\mu}_{mn}f_{mn}+\hat{\nu}
	\end{equation} 
which is a convex constraint since $\lambda_{mn}<0$. Again, if the inequality in \eqref{eq:pumprel} is satisfied with equality for all pumps, the relaxation of \eqref{eq:pumphead} to \eqref{eq:pumprel} is deemed as exact.

%Interestingly, the relaxation of \eqref{eq:headloss}--\eqref{eq:pumphead} to \eqref{eq:relaxed}--\eqref{eq:pumprel} was proved to be exact for the OWF problem in \cite{TaylorOWFCDC}. This exactness can be attributed to the particular cost function of the OWF problem as well as four simplifying assumptions: \emph{a1)} no pumps on cycles; \emph{a2)} known flow directions; \emph{a3)} if a node is fed by more than one pipes, all of them have to be equipped with valves; and \emph{a4)} no fixed-speed pumps. Assumption \emph{a2)} could be possibly met if the WDS operator has a good sense of the network conditions, whereas \emph{a3)} may be practically restricting. If one tries to tackle \eqref{eq:W1} as an instance of the OWF problem of \cite{TaylorOWFCDC}, the relaxation will not be exact. This is because even if a \eqref{eq:W1} instance satisfies \emph{a1)--a3)}, assumption \emph{a4)} is violated since pump speeds are fixed parameters for \eqref{eq:W1}. Therefore, problem \eqref{eq:W1} violates the sufficient conditions of \cite{TaylorOWFCDC} for the relaxation to be exact. In addition, numerical tests on solving WF instances even in simple water networks demonstrate that the approach of \cite{TaylorOWFCDC} fails to solve \eqref{eq:W1}, i.e., the convex relaxation is not exact. Because of this, we next pursue an alternative approach to render the relaxation exact when solving the WF problem. The suggested approach succeeds in a broader class of WDS configurations compared to those considered in \cite{TaylorOWFCDC} and \cite{TaylorOWFcones}. 

One may now try solving the feasibility problem in \eqref{eq:W1} by replacing \eqref{eq:headloss}--\eqref{eq:pumphead} by \eqref{eq:relaxed}--\eqref{eq:pumprel}. If all the relaxed constraints are satisfied with equality, the relaxation is \emph{exact}. Unfortunately, toy numerical tests can verify that such relaxation is inexact even for simple water networks. To favor exact relaxations, we convert the feasibility to a minimization problem with a judiciously selected cost. We pose the MI-QCQP
	\begin{align*}\label{eq:W2}
	\min_{\bef,\bh,\bx}~&~s(\bh)\tag{W2}\\
%	\mathrm{over}~& ~ \bef,\bh,\bx\\
	\mathrm{s.to}~&~ \eqref{eq:wmc},\eqref{eq:relaxed},\eqref{eq:pumprel}
	\end{align*}
where the binary vector $\bx$ contains all $\{x_{mn}\}_{(m,n)\in\bmcP_a}$, and the cost is defined as
	\begin{equation*}%\label{eq:gh}
	s(\bh):=\sum_{(m,n)\in \bmcP_a}|h_m-h_n|-\sum_{(m,n)\in \mcP_a}(h_n-h_m).
	\end{equation*}
The cost	sums up the absolute pressure drops across lossy pipes minus the pressure gains added by pumps. This is in contrast to the MI-QCQP of \cite{Singh18ACC} that recovers a gas flow solution, whose penalty did not include compressors. It will be shown in Section~\ref{subsec:exact} that the MI-QCQP formulation \eqref{eq:W2} guarantees exactness under realistic network conditions.

Problem \eqref{eq:W2} is non-convex due to the binary variables $\bx$. Given the advancements in MI-QCQP solvers, this minimization can be handled for moderately sized networks. Recall that every lossy pipe is associated with one binary and one continuous optimization variable. To accelerate the computations, we next provide a simple pre-processing step to determine the flows in all edges not belonging to a cycle.

\begin{lemma}\label{le:fixflow}
Let $\bef$ be the unique flow solution of \eqref{eq:W1}, and consider any vector $\tbf$ satisfying $\bA^\top\tbf=\bd$. For any edge $(m,n)$ not belonging to a cycle, it holds that $\tilde{f}_{mn}=f_{mn}$.
\end{lemma}

\begin{IEEEproof}
Consider the minimum-norm solution $\bef_0:=\left(\bA^\top\right)^\dagger\bd$ to the linear system of \eqref{eq:wmc}, where $\left(\bA^\top\right)^\dagger$ is the pseudo-inverse of $\bA^\top$. Any other solution of \eqref{eq:wmc} can be expressed as $\bef=\bef_0+\bn$ for some vector $\bn\in\nullspace(\bA^\top)$. 

The space $\nullspace(\bA^\top)$ can be represented using the cycles of graph $\mcG$ as explained next. For cycle $\mcC$ in $\mcG$, select an arbitrary direction and define its indicator vector $\bn_\mcC\in\mathbb{R}^{P}$ as
\begin{equation*}
n_{\mcC,p}:=
\begin{cases}
	0&,~\text{if edge } p\notin \mcC\\
	+1&,~\text{if the directions of }p\text{ and }\mcC\text{ coincide}\\
	-1&,~\text{otherwise.}\\
\end{cases}
\end{equation*}
We are particularly interested in a set of \emph{fundamental cycles} defined as follows~\cite{GodsilRoyle}: In graph $\mcG=(\mcN,\mcP)$, select a spanning tree $\mcT$. Every edge $p\in\mcP\setminus \mcT$ along with some edges of $\mcT$ form a cycle. The cycles formed using all edges in $\mcP\setminus \mcT$ comprise the set $\mcL_\mcT$ of fundamental cycles. In general, a graph may have multiple spanning trees, making the set of fundamental cycles non-unique.

A key property of fundamental cycles is that the space $\nullspace(\bA^\top)$ is spanned by the indicator vectors for any set of fundamental cycles~\cite[Corollary~14.2.3]{GodsilRoyle}. Therefore, the vector $\bn$ in $\bef=\bef_0+\bn$ can be decomposed as
\begin{equation}\label{eq:basis}
\bn=\sum_{\ell \in \mcL_\mcT}\alpha_\ell\bn_\ell
\end{equation}
where $\bn_\ell$ is the indicator vector for cycle $\ell$ and $\alpha_\ell\in\mathbb{R}~\forall \ell$. 

Consider the $p$-th entry of $\bn$. If edge $p\in\mcP$ does not belong to any cycle, then it does not belong to any fundamental cycle. Then $n_{\ell,p}=0$ for all $\ell\in\mcL_\mcT$, and $n_p=0$ follows from \eqref{eq:basis}. To conclude, if edge $p$ does not belong to a cycle, then all solutions $\bef$ of \eqref{eq:wmc} agree in their $p$-th entry.
\end{IEEEproof}

Lemma~\ref{le:fixflow} ensures that for weakly meshed WDS, the number of variables in \eqref{eq:W2} can be reduced markedly. Moreover, the flow on any edge not belonging to a cycle coincides with the related entry of the minimum-norm solution to \eqref{eq:wmc}. 

%%%%%%%%%%%%%%%%%%%%%%%%%%%%%%%%%%%%%%%%%%%
\subsection{Exactness of the Relaxation}\label{subsec:exact}
The next result provides conditions under which a minimizer of \eqref{eq:W2} satisfies \eqref{eq:relaxed}--\eqref{eq:pumprel} with equality.

\begin{theorem}\label{th:wf1}
In a WDS where no edge $(m,n)\in\mcP$ belongs to more than one cycle, a minimizer of \eqref{eq:W2} minimizes \eqref{eq:W1} as well, if \eqref{eq:W1} is feasible.
\end{theorem}

\begin{figure}[t]
	\centering
	\includegraphics[scale=0.7]{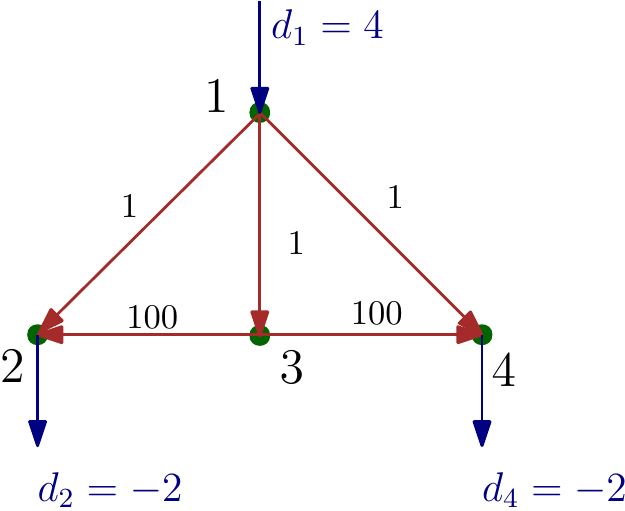}
	\caption{A pathological WDS for which the relaxation of \eqref{eq:W2} is inexact.}
	\label{fig:fail}
\end{figure}

Theorem~\ref{th:wf1} (shown in the appendix) asserts that the MI-QCQP relaxation of \eqref{eq:W1} to \eqref{eq:W2} is exact in water networks with non-overlapping cycles, that is cycles sharing no edges. The claim holds regardless of water demand or the presence of pumps. This is in contrast to the exactness claims of \cite{Singh18ACC}, where compressors were not allowed on cycles. Although the assumptions of Theorem~\ref{th:wf1} may not always hold, the numerical tests of Section~\ref{sec:tests} indicate the relaxation is exact and \eqref{eq:W2} solves \eqref{eq:W1} for practical WDS with overlapping cycles.

The condition of Theorem~\ref{th:wf1} cannot be relaxed analytically: One could construct pathological WDS with overlapping cycles that render the relaxation of \eqref{eq:W2} inexact. To present such a counterexample, consider the $4$-node and $5$-pipe WDS of Fig.~\ref{fig:fail}. The coefficients $c_{mn}$'s are shown on the respective pipes. Nodes $2$ and $4$ host water demands, and the reference node $1$ supplies water at the reference pressure of $h_1=10$. The edge $(1,3)$ belongs to two cycles and hence the condition of Theorem~\ref{th:wf1} is violated. The minimizer of \eqref{eq:W2} satisfies the relaxed constraint \eqref{eq:relaxed} with strict inequality. The pressures at nodes $1$ and $3$ were found to differ by $h_1-h_3=3.435$, while the frictional drop was $c_{13}f_{13}^2=0.014$. This WDS though is characterized by a strong disparity of pipe coefficients. For a WDS with limited variations in pipe dimensions and material, such disparity is not anticipated in practice. Similar to the example of Fig.~\ref{fig:fail}, one can construct setups with multiple fixed-pressure nodes for which the MI-QCQP relaxation is inexact.
% \cmr{[Some $\frac{\max c_p}{\min c_p}$ ratios from EPANET to support this?]}\cmb{EPANET numbers do not help here} 

%%Therefore, the numerical tests conducted on practical networks with overlapping cycles in Section~\ref{sec:tests} yield exact relaxation with \eqref{eq:W2}.
	
Since the WDS of Fig.~\ref{fig:fail} does not host pumps, the solution to the previous WF problem was eventually found via \eqref{eq:cWF}. This motivates us to exploit the ability of the energy function-based approach to handle overlapping cycles alongside the merit of the MI-QCQP relaxation to handle pumps in cycles. To solve the WF problem for a broader class of WDS network topologies, a hybrid solver is deferred to Section~\ref{sec:hybrid}. Before that, the formulation of \eqref{eq:W2} is contrasted to prior work.

\subsection{Comparison to Prior Work and Extensions}\label{subsec:comparison}
The proposed WF solver involves three ingredients:
\renewcommand{\labelenumi}{\emph{I\arabic{enumi})}}
\begin{enumerate}
\item replacing \eqref{eq:headloss} by the MIQP model of \eqref{eq:relaxed};
\item replacing \eqref{eq:pumphead0} by the convex constraint of \eqref{eq:pumprel}; and
\item solving \eqref{eq:W2} in lieu of \eqref{eq:W1}. 
\end{enumerate}

The key contribution of this section is \emph{I3)}, since \emph{I1)} and \emph{I2)} have appeared before: References \cite{TaylorOWFCDC}--\cite{TaylorOWFcones} and \cite{ZamzamOWPF} consider the OWF problem assuming known flow directions (say $f_{mn}\geq 0$ for all pipes), and relax \eqref{eq:headloss} to [cf.~\eqref{eq:relaxed}]
\begin{equation}\label{eq:taylorpipe}
h_m-h_n\geq c_{mn}f_{mn}^2,~f_{mn}\geq 0.
\end{equation}
Regarding pumps, references \cite{TaylorOWFCDC}--\cite{TaylorOWFcones} and \cite{ZamzamOWPF} consider only variable-speed pumps. Moreover, each pump speed $\omega_{mn}$ is limited within $0\leq \omega_{mn}\leq \bar{\omega}_{mn}$, which may be unrealistic since pumps come with positive lower speed bounds. Substituting this speed range into the monotonic formula of \eqref{eq:pumphead} yields the convex quadratic constraint
\begin{align}\label{eq:taylorpump}
h_n-h_m\leq\lambda_{mn} f_{mn}^2+\mu_{mn}\bar{\omega}_{mn} f_{mn} + \nu_{mn}\bar{\omega}^2_{mn}.
\end{align}
Heed that the actual speed $\omega_{mn}$ has been eliminated from \eqref{eq:taylorpump}. Nevertheless, once an OWF minimizer is found [and so a triplet $(h_m,h_n,f_{mn})$ satisfying \eqref{eq:taylorpump}], the speed $\omega_{mn}$ making \eqref{eq:taylorpump} an equality can be readily recovered~\cite{TaylorOWFCDC}. 

References \cite{TaylorOWFCDC}--\cite{TaylorOWFcones} aim at minimizing the energy losses across pipes. Ignoring the details of geographical elevation and electricity prices, their cost simplifies to $\sum_{(m,n)\in \bmcP_a}c_{mn}f_{mn}^3$. Thanks to its form, minimizing this cost subject to \eqref{eq:taylorpipe}--\eqref{eq:taylorpump} renders the relaxation in \eqref{eq:taylorpipe} numerically exact. This was also analytically shown granted the simplifying assumptions of: \emph{a1)} no pumps on cycles; \emph{a2)} known flow directions; \emph{a3)} if a node is fed by more than one pipes, all of them have to be equipped with valves; and \emph{a4)} no fixed-speed pumps. Reference~\cite{ZamzamOWPF} considers the joint scheduling of WDS and electric power distribution networks. Since the exactness of \eqref{eq:taylorpipe} is not promoted by the cost anymore, a feasible point pursuit-based approach is applied to find a stationary point satisfying \eqref{eq:taylorpipe} with equality.

From the preceding review, it is evident that albeit a relaxation of \eqref{eq:headloss} is convenient, its exactness is not always guaranteed. While some OWF instances feature exactness \cite{TaylorOWFCDC}--\cite{TaylorOWFcones}; others call for additional measures~\cite{ZamzamOWPF}. We shall next address the natural question: \emph{Q1)} Can one use one of the aforesaid OWF solvers to solve the WF task? Numerical tests on WF instances even in simple WDS demonstrate that the approaches of \cite{TaylorOWFCDC}--\cite{TaylorOWFcones} and \cite{ZamzamOWPF} fail to solve \eqref{eq:W1}, that is the convex relaxation is not exact. This is because even if an instance of \eqref{eq:W1} satisfies all other conditions needed for each one of these solvers, the common assumption \emph{a4)} is violated since pump speeds are fixed for \eqref{eq:W1}.

One may reverse question \emph{Q1)} and pose \emph{Q2)}: Can the developed WF solvers be used towards solving OWF tasks? One could identify two possible uses. First, the proposed tools can be used instead of EPANET in existing zero-order OWF algorithms that rely on a WF solver. Second, the success of the proposed penalized relaxation could be adopted in OWF tasks. For example, our previous work~\cite{singh2018optimal} deals with optimal pump scheduling under dynamic electricity pricing. This OWF task handles unknown flow directions and is relaxed to an MI-SOCP after $g_{mn}(f_{mn};\omega^0_{mn})$ in \eqref{eq:pumphead} is approximated as constant. The relaxation is provably exact under certain conditions and after adding a penalty to the cost of the MI-SOCP. The relaxation is numerically exact under a broader range of WDS conditions. Interestingly, the penalty that worked for \cite{singh2018optimal} does not work for the WF task studied here. For this reason, this work put forth the penalty $s(\bh)$ in \eqref{eq:W2} and took a totally different route for establishing exactness. The sufficient conditions for exactness of Theorem~\ref{th:wf1} are significantly simpler compared to those in prior works on OWF. 

% assuming fixed-speed pumps and small $\lambda_{mn}$ and $\mu_{mn}$

%%%%%%%%%%%%%%%%%%%%%%%%%%%%%%%%%%%%%%%%%%%%%%%%%%%
\begin{figure}[t]
	\centering
	\includegraphics[scale=0.6]{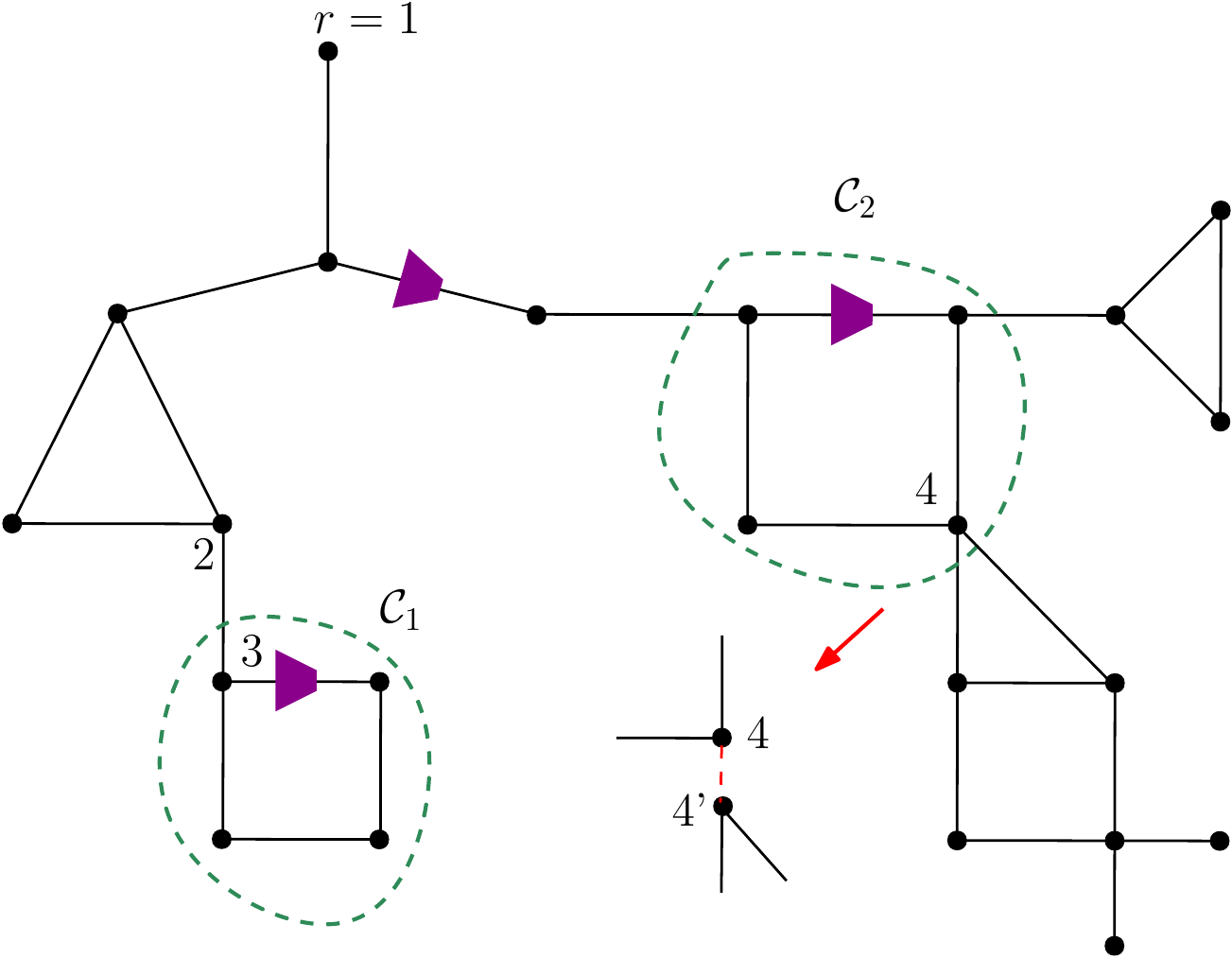}\vspace*{1em}
	\includegraphics[scale=0.6]{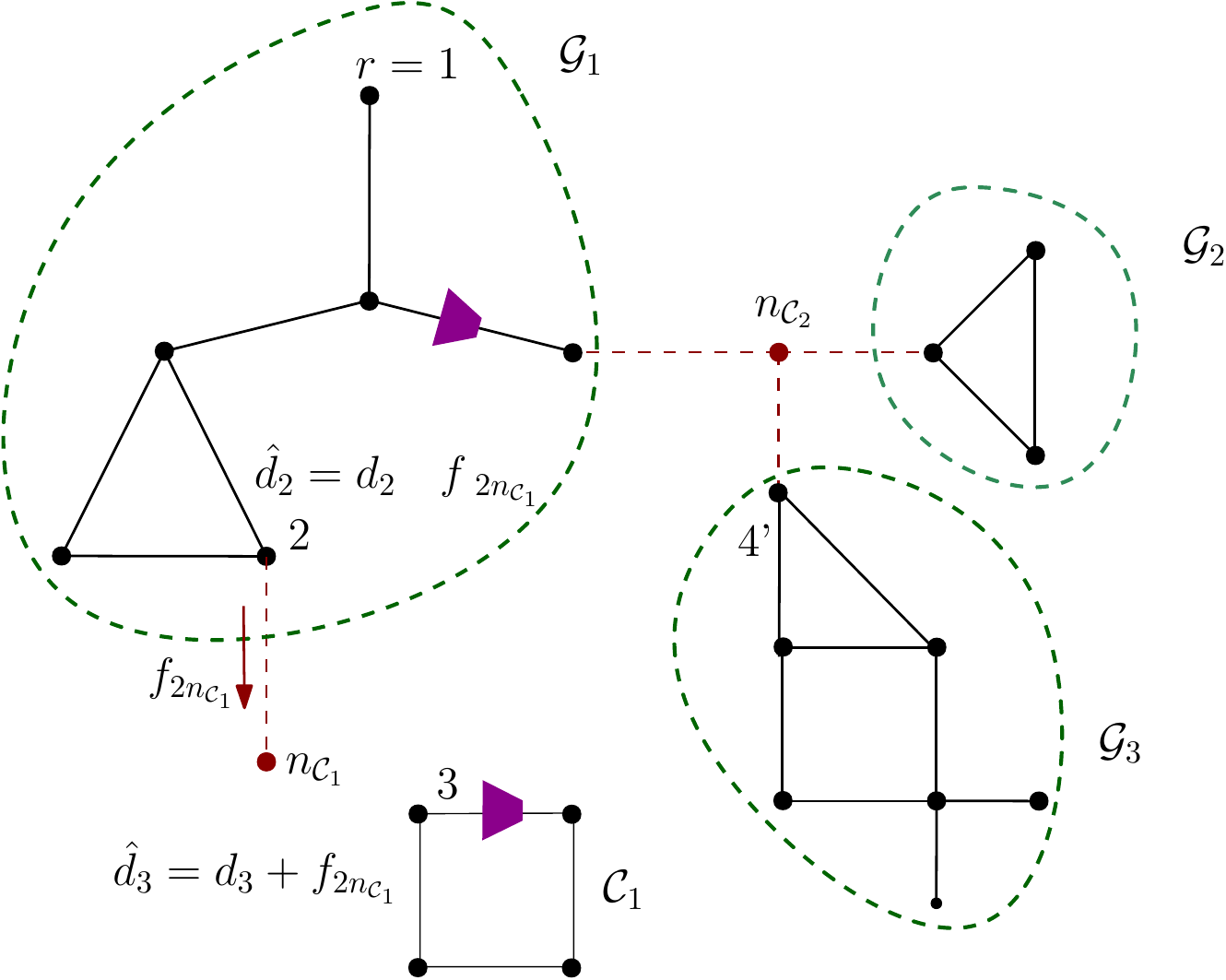}
	\caption{\emph{Top:} A WDS with pumps in non-overlapping cycles; \emph{Bottom}: Connected components after removing cycles $\mcC_1$ and $\mcC_2$ that carry pumps.}
	\label{fig:hybrid}
\end{figure}

\section{A Hybrid Water Flow Solver}\label{sec:hybrid}
The hybrid WF solver relaxes the assumption of non-overlapping cycles of Theorem~\ref{th:wf1} to the following condition.

\begin{assumption}\label{as:A1}
The water network has no pumps in overlapping cycles.
\end{assumption}

The assumption permits overlapping cycles, but these particular cycles should carry no pumps. For a network satisfying Assumption~\ref{as:A1}, the WF task can be solved through the following steps illustrated also in Fig.~\ref{fig:hybrid}:
		\renewcommand{\labelenumi}{\emph{T\arabic{enumi})}}
	\begin{enumerate}
		\item Any cycle $\mcC$ with pumps is non-overlapping. Thus, for any node $n$ belonging to $\mcC$, two cases arise:\\
		\emph{i)} Node $n$ belongs only to cycle $\mcC$ (node $3$ of $\mcC_1$); or\\
		\emph{ii)} Node $n$ belongs to other cycle(s) $\mcC'$ as well; yet $\mcC$ and $\mcC'$ share no edge (e.g., node $4$ of $\mcC_2$).\\
		Then, reduce graph $\mcG$ to $\mcG'$ through the steps:
		\begin{itemize}
			\item  If for a cycle $\mcC$ having pumps, all nodes do not belong to any other cycle, replace $\mcC$ by a supernode $n_\mcC$. 
			\item If for a cycle $\mcC$ having pumps, node $n$ belongs to other cycle(s), identify the edges $(n_1,n)$ and $(n,n_2)$, which belong to $\mcC$. These two edges belong to no other cycles as $\mcC$ is non-overlapping. Split the node $n$ into $n$ and $n'$ connected by a lossless edge $(n,n')$, such that all edges in $\mcG$ other than $(n_1,n)$ and $(n,n_2)$ that were incident on $n$ are now incident on $n'$. After repeating this step for all nodes in $\mcC$ that belong to multiple cycles, replace $\mcC$ by a supernode $n_\mcC$. 
		\end{itemize}
	This  process ensures that in $\mcG'$ all supernodes and the pumps left out from $\mcG$ do not appear on cycles.
		\item Use Lemma~\ref{le:fixflow} to compute the water flows on the edges $(m,n_\mcC)$ incident to supernodes $n_\mcC$'s in $\mcG'$.  
		\item For each edge $(m,n_\mcC)$ in $\mcG'$, modify the injection at node $m$ as $\hat{d}_m=d_m-f_{mn_\mcC}$.
		\item Partition $\mcG'$ into a set of connected components $\mcG_c$ by removing the supernodes and their incident edges. Each $\mcG_c$ has known injections and bears no pumps in cycles. 
		\item Use the stitching algorithm of Section~\ref{subsec:stitching} per component $\mcG_c$ to find the flows within $\mcG_c$. 
		\item Each supernode $n_\mcC$ is split back to the cycle $\mcC$ it replaced. If the edge $(m,n_\mcC)$ of $\mcG'$ corresponded to edge $(m,n)$ of $\mcG$, modify the injection at node $n\in\mcC$ as $\hat{d}_n=d_n+f_{mn_\mcC}$.
		\item Solve \eqref{eq:W2} per cycle $\mcC$ to find the water flows on $\mcC$.
		\item Given vector $\bef$, find the pressure vector $\bh$ using Lemma~\ref{le:fh}.
	\end{enumerate}

%\begin{figure}[t]
%	\centering
%	%\hspace*{-3em}
%	\includegraphics[scale=0.6]{Example3_1.eps}\vspace*{1em}
%	%\hspace*{-3em}
%	\includegraphics[scale=0.6]{Example3_2.eps}
%	%\vspace*{-1em}
%	\caption{\emph{Top:} A WDS with pumps in non-overlapping cycles; \emph{Bottom}: Connected components after removing cycles $\mcC_1$ and $\mcC_2$ that carry pumps.}
%	\label{fig:hybrid}
%\end{figure}

Let us apply the previous steps on the $23$-node and $28$-pipe water network of Fig.~\ref{fig:hybrid}. There are two cycles carrying pumps, marked as $\mcC_1$ and $\mcC_2$. Node $4$ of cycle $\mcC_1$ belongs to multiple cycles and hence it is split in $4$ and $4'$. Next, removing $\mcC_1$ and $\mcC_2$ along with the edges that connect these cycles with the rest of the graph results in the connected subgraphs $\mcG_1$, $\mcG_2$, and $\mcG_3$ shown on the bottom of Fig.~\ref{fig:hybrid}. The demands on boundary nodes are modified as per steps \emph{T3)} and \emph{T6)}. The edge flows within each connected component are subsequently found by solving \eqref{eq:cWF}. The flows on $\mcC_1$ and $\mcC_2$ are finally found by \eqref{eq:W2}. 

This WDS setup could not be handled by the energy function-based approach of Section~\ref{sec:energy} alone due to the presence of pumps on cycles $\mcC_1$ and $\mcC_2$. The convex relaxation of Section~\ref{sec:MI-QCQP} alone is not guaranteed to succeed either, due to the presence of overlapping cycles. Combining the merits of each method and leveraging Lemma~\ref{le:fixflow}, this hybrid method can handle successfully this WDS setup. 

%%%%%%%%%%%%%%%%%%%%%%%%%%%%%%%%%%%%%%%%%%%%%%%%%%%%
\section{Numerical Tests}\label{sec:tests} 

\begin{figure}[t]
	\centering
	\includegraphics[scale=0.46]{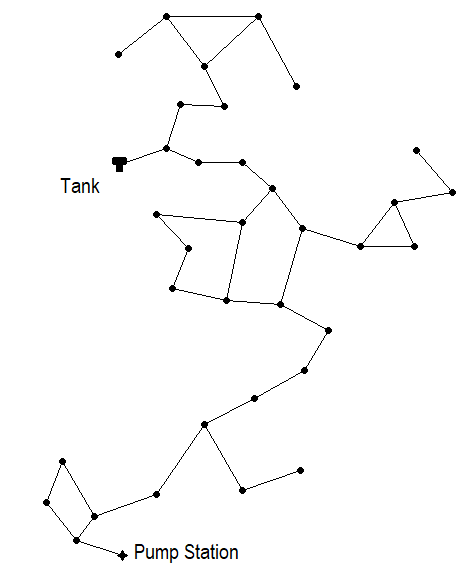}
	\caption{EPANET Example Network-2 of a WDS from Cherry Hills, CT~\cite{Rossman94net2}.}
	\label{fig:net2map}
\end{figure}

The new WF solvers were evaluated under the different network configurations of Fig.~\ref{fig:flowchart}. First, the new WF solvers were evaluated on the EPANET \emph{Example Network-2} representing a water network from Cherry Hills, Connecticut shown in Fig.~\ref{fig:net2map}~\cite{Rossman94net2}. This network consists of $P=40$ pipes, $N=34$ demand nodes, one tank, and one pump station. We modified the network by representing the pump station as a reservoir with pressure $100$~ft connected to a fixed-speed pump. All nodes were assumed at the same elevation. The pipe friction coefficients $c_{mn}$'s, and base demand vector $\bd$, were derived from the EPANET benchmark. Note that the WDS in Fig.~\ref{fig:net2map} has overlapping cycles but the pump is not in a loop. Thus, it qualifies for being solved using the energy function method alongside the stitching algorithm of Section~\ref{subsec:stitching}.

\begin{figure}[t]
	\centering
	\includegraphics[scale=0.19]{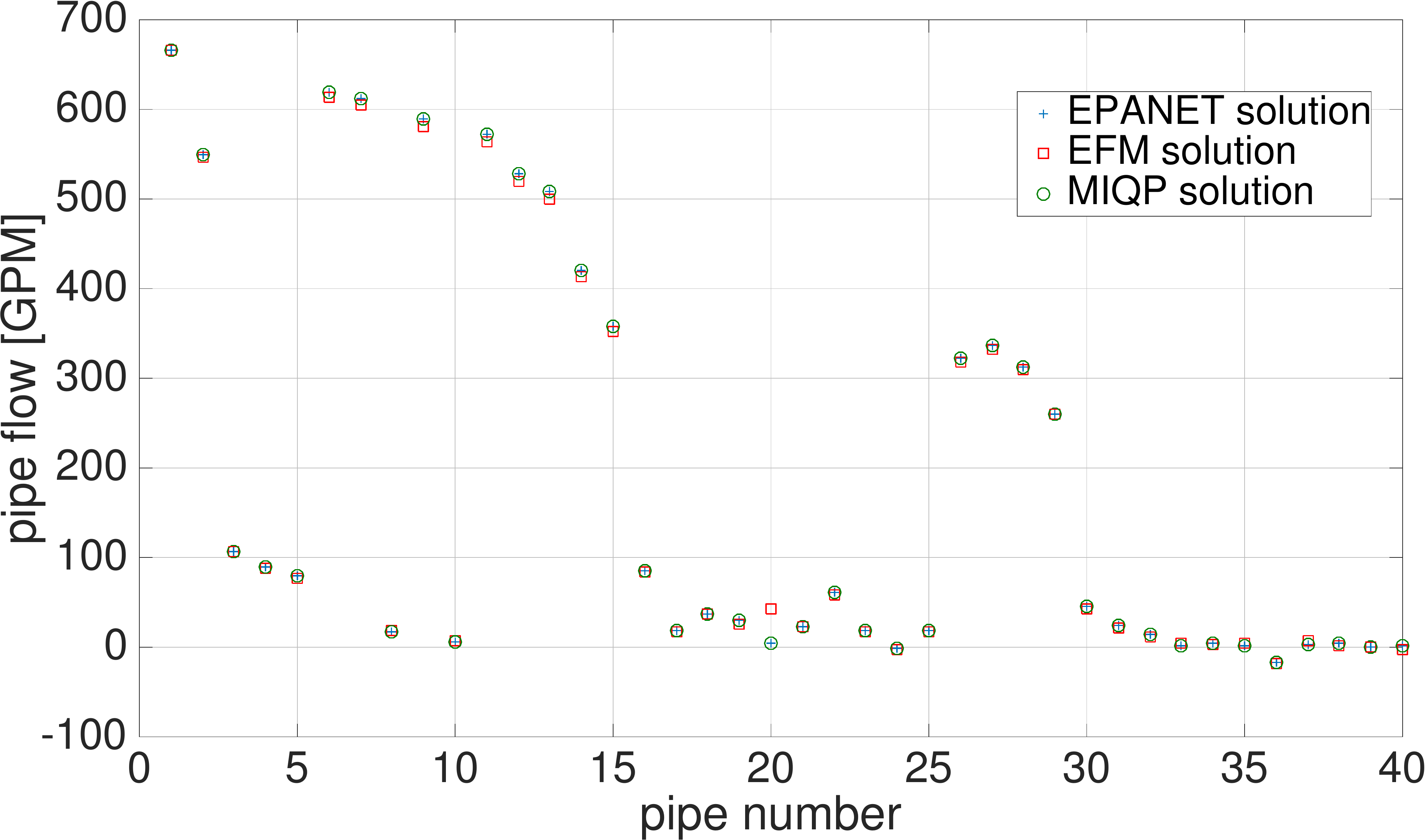}
	\caption{Water flows obtained by the EPANET solver, the energy function minimization of~\eqref{eq:cWF}, and the MI-QCQP of~\eqref{eq:W2} for the WDS of Fig.~\ref{fig:net2map}.}
	\label{fig:cr_epa_ef}
\end{figure}

We tested whether the proposed solvers yield the same solution as EPANET. For this purpose, we obtained WF solutions using the constrained energy function minimization of \eqref{eq:cWF} and the MI-QCQP of \eqref{eq:W2}. Problem \eqref{eq:cWF} was solved using the closed-form dual decomposition steps of Section~\ref{subsec:cWF} for $\mu=10^{-4}$ within $20,000$ iterations. The MI-QCQP in \eqref{eq:W2} was solved using the MATLAB-based toolbox YALMIP along with the mixed-integer solver CPLEX~\cite{yalmip}, \cite{cplex}. All tests were run on a $2.7$~GHz, Intel Core i5 computer with $8$~GB RAM. The flows obtained by the two solvers were very close to the EPANET solution, as illustrated in Fig.~\ref{fig:cr_epa_ef}. EPANET uses more detailed flow models, e.g., the coefficient $c_{mn}$ in \eqref{eq:headloss} depends weekly on flow $f_{mn}$. The $20,000$ iterations for the primal-dual updates of \eqref{eq:cWF} were completed within $4.5$~sec, and the running time of the MI-QCQP~\eqref{eq:W2} was $1.6$~sec.

\begin{figure}[t]
	\centering
	%\hspace*{-3em}
	\includegraphics[scale=0.2]{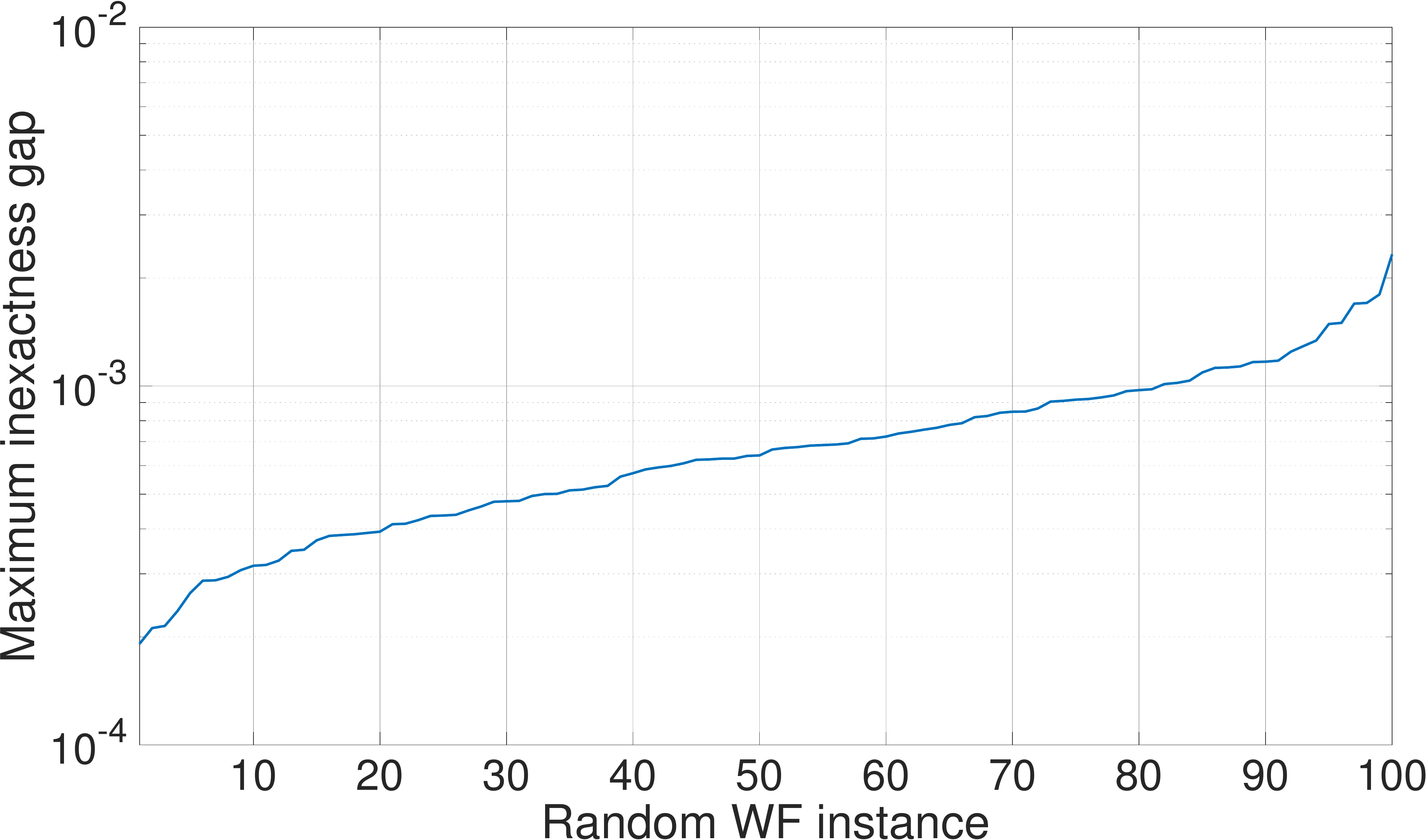}
	%\hspace*{-3em}
	\includegraphics[scale=0.2]{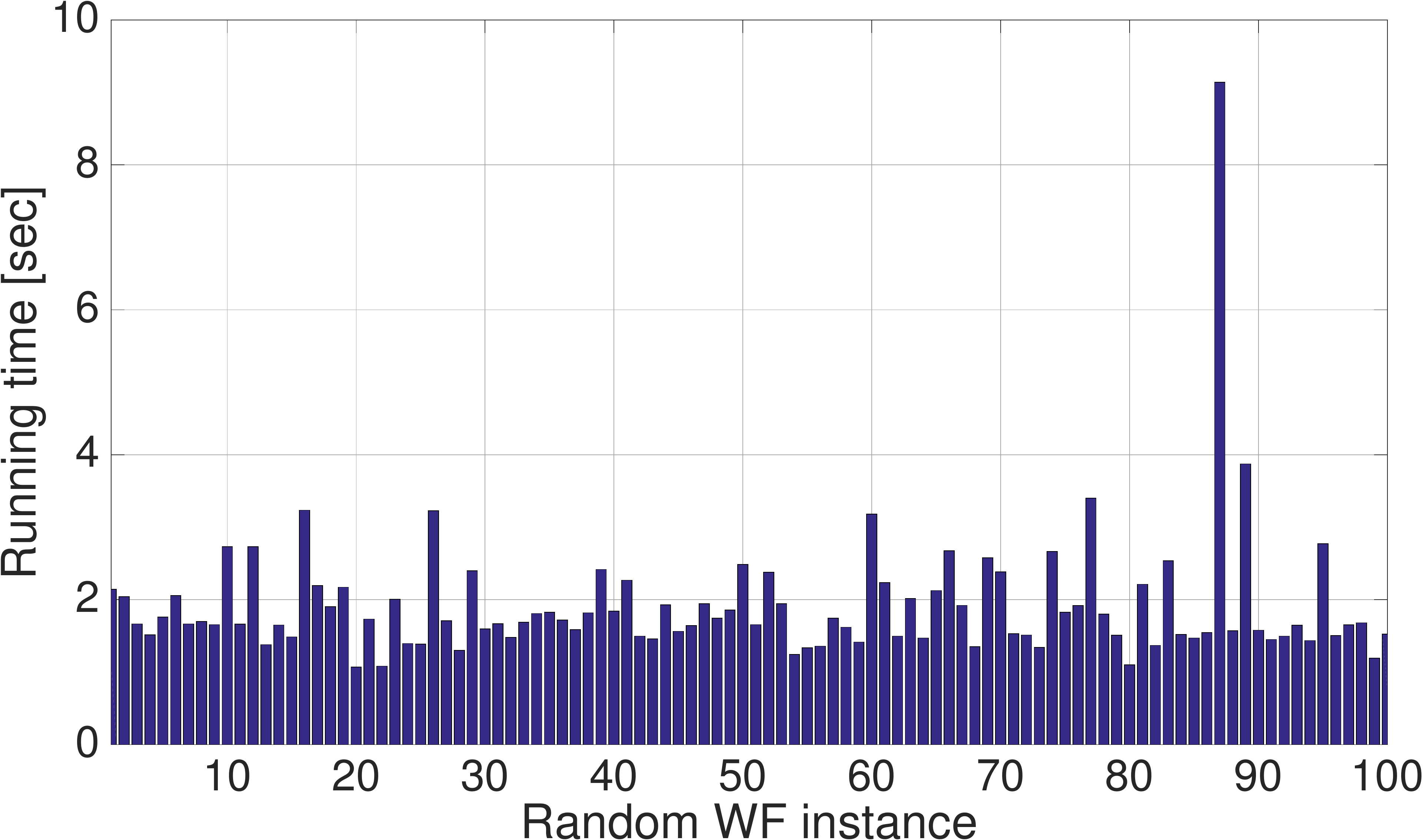}
	%\vspace*{-1em}
	\caption{\emph{Top:} Maximum inexactness gap attained by \eqref{eq:W2} over $100$ random WF instances. \emph{Bottom}: Running time for \eqref{eq:W2} over random WF instances.}
	\label{fig:CR_MC}
\end{figure}

We next evaluated the performance of the MI-QCQP-based solver of \eqref{eq:W2} in finding the correct WF solution, when the conditions of Theorem~\ref{th:wf1} are not satisfied. Specifically, the network of Fig.~\ref{fig:net2map} has overlapping cycles. To represent various demand levels, we generated $100$ random WF instances by scaling the benchmark demand $\bd$ by a scalar uniformly drawn from $[0,1.5]$. Given a minimizer of \eqref{eq:W2}, we defined the inexactness over lossy pipe $(m,n)$ as ${|h_m-h_n|-c_{mn}f_{mn}^2}$. For each run of \eqref{eq:W2} with a random input, the ranked maximum inexactness gap over all lossy pipes is displayed on Fig.~\ref{fig:CR_MC} (top). Despite violating the conditions of Th.~\ref{th:wf1}, the inexactness gap was small for all random tests. Computationally, the running time of \eqref{eq:W2} over the $100$ instances with a median value of only $1.68$~sec; see Fig.~\ref{fig:CR_MC} (bottom).

\begin{figure}[t]
	\centering
	\includegraphics[scale=0.3]{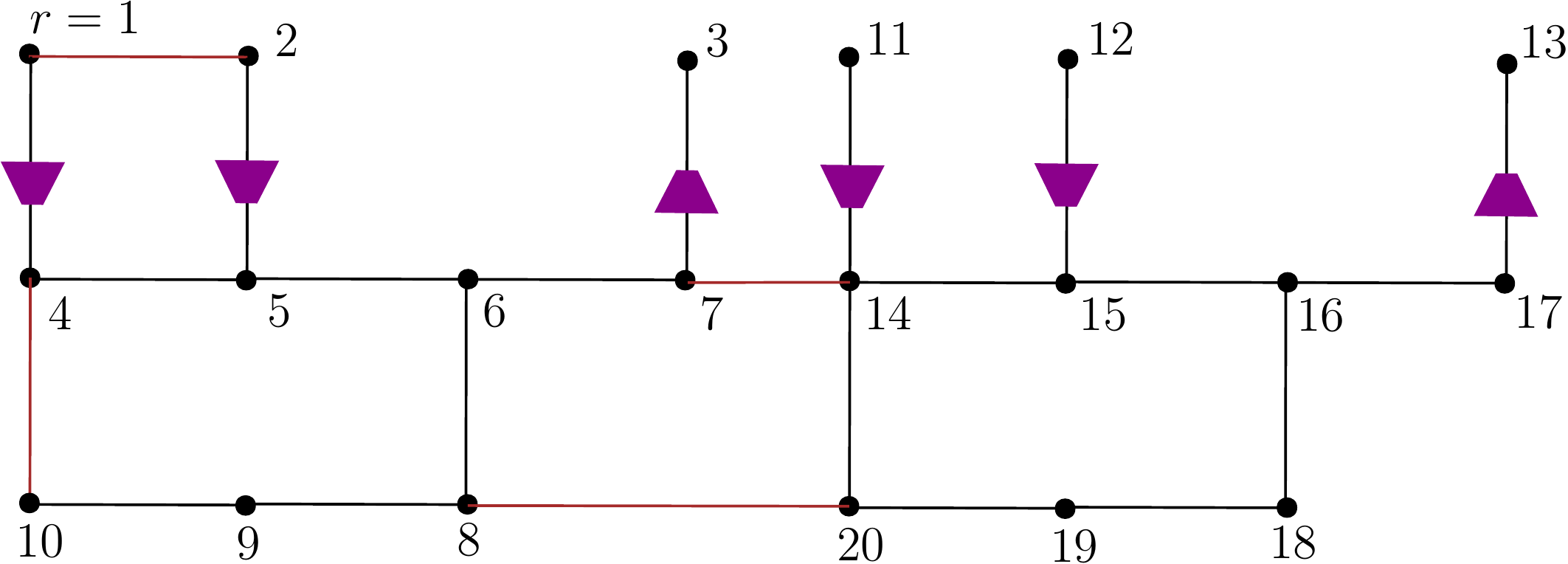}
	\caption{Benchmark 20-node WDS derived from the WDS of \cite{singh2018optimal}.}
	\label{fig:WDS20}
	%	\vspace*{-1em}
\end{figure}

Our WF solvers were subsequently tested on the 20-node synthetic network of Fig.~\ref{fig:WDS20} and under different demands. The WDS of Fig.~\ref{fig:WDS20} was created by combining two copies of the popular 10-node WDS representing the modified Arava Valley network from Israel; see \cite{singh2018optimal}, \cite{CohenQHmodel} for network data. The additional pipes $(1,2)$, $(7,14)$, and $(8,20)$ have the same dimensions as pipes $(4,5)$, $(5,6)$, and $(6,8)$, respectively. These were added to generate different network conditions. The parameters for all the pumps were derived by scaling the pump parameters of \cite{TaylorOWFcones} by $0.25$, and were set to $(\lambda_{mn}, \bar{\mu}_{mn}, \bar{\nu}_{mn})=(-2.735\cdot 10^{-5}, 0.0129, 55.83)$. 

The performance of \eqref{eq:W2} was tested for two network configurations derived from the WDS of Fig.~\ref{fig:WDS20}: \emph{C-i)} pipe $(4,10)$ removed; \emph{C-ii)} pipe $(8,20)$ removed. Note that configuration \emph{C-i)} has pumps in non-overlapping cycle along side other overlapping cycles not hosting pumps. Thus, the WF problem for \emph{C-i)} can provably be solved using the hybrid WF solver of Section~\ref{sec:hybrid}. Moreover, configuration \emph{C-ii)} has pumps in overlapping cycles and thus cannot be provably solved even by the hybrid WF solver. However, we numerically tried \eqref{eq:W2} in solving $100$ random WF instances for both \emph{C-i)} and \emph{C-ii)}.

The $100$ random WF instances were generated by fixing pressure $h_1=10$ m and drawing the remaining pressures independently as Gaussian random variables of mean $10$~m and variance $2$ m$^2$, or $\mcN(10,2)$. The on/off statuses of pumps were drawn as independent Bernoulli random variables  with mean of $0.5$. The flow for on pumps was drawn with uniform probability on the allowable pump flow $[250,1500]$ m$^3/$hr. For off pumps, flows were drawn from $\mcN(0,200)$. The pressure added by pumps was then calculated using \eqref{eq:pumphead}. For pump $(m,n)$, the receiving node pressure was updated as $h_n:=h_m+g_{mn}$. The water flow in all lossy pipes was calculated from the so obtained pressures and \eqref{eq:headloss}. Once the complete flow vector was obtained, the injections were computed from \eqref{eq:1}. The obtained vector of injections, pump on/off status, and pressure $h_1$ served as a feasible input for the WF task. Similarly $100$ feasible random WF instances were generated for configuration \emph{C-ii)}, and solved using \eqref{eq:W2}. The maximum time for solving \eqref{eq:W2} was set to $1$ min in CPLEX. 

\begin{figure}[t]
	\centering
	\includegraphics[scale=0.2]{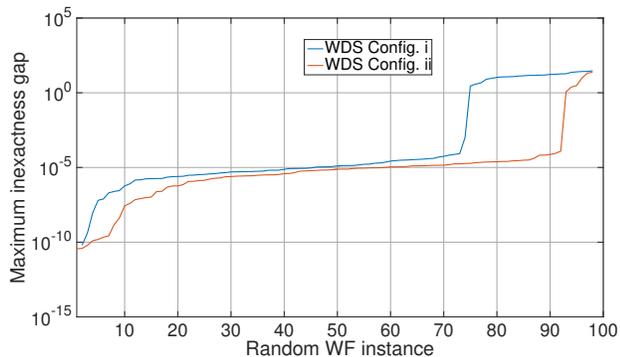}
	\caption{Maximum inexactness gap attained by \eqref{eq:W2} over 100 random WF instances for configurations \emph{C-i)} and \emph{C-ii)} of the WDS of Fig.~\ref{fig:WDS20}.}
	\label{fig:wds20_MC}
	%	\vspace*{-1em}
\end{figure}

The WF task for both configurations was solved within a $1$-minute deadline for $98$ out of the $100$ random instances with the average time being $1.01$ and $1.06$~sec for \emph{C-i)} and \emph{C-ii)}, respectively. The ranked maximum inexactness gap over all lossy pipes for \emph{C-i)} and \emph{C-ii)} is shown in Fig.~\ref{fig:wds20_MC}. Although the conditions of Th.~\ref{th:wf1} were not satisfied by either configuration, the maximum inexactness obtained was less than $10^{-3}$ for $74$ and $92$ out of the $98$ instances of \emph{C-i)} and \emph{C-ii} that were solved. For the two cases where the MI-QCQP failed to converge within the $1$-min deadline, we modified the big-$M$ parameter value from $M=300$ to $M=80$. The result was that the two instances were solved in $0.13$ and $0.82$ seconds with the inexactness gap being less than $10^{-5}$.

\begin{figure}[t]
	\centering
	%\hspace*{-3em}
	\includegraphics[scale=0.2]{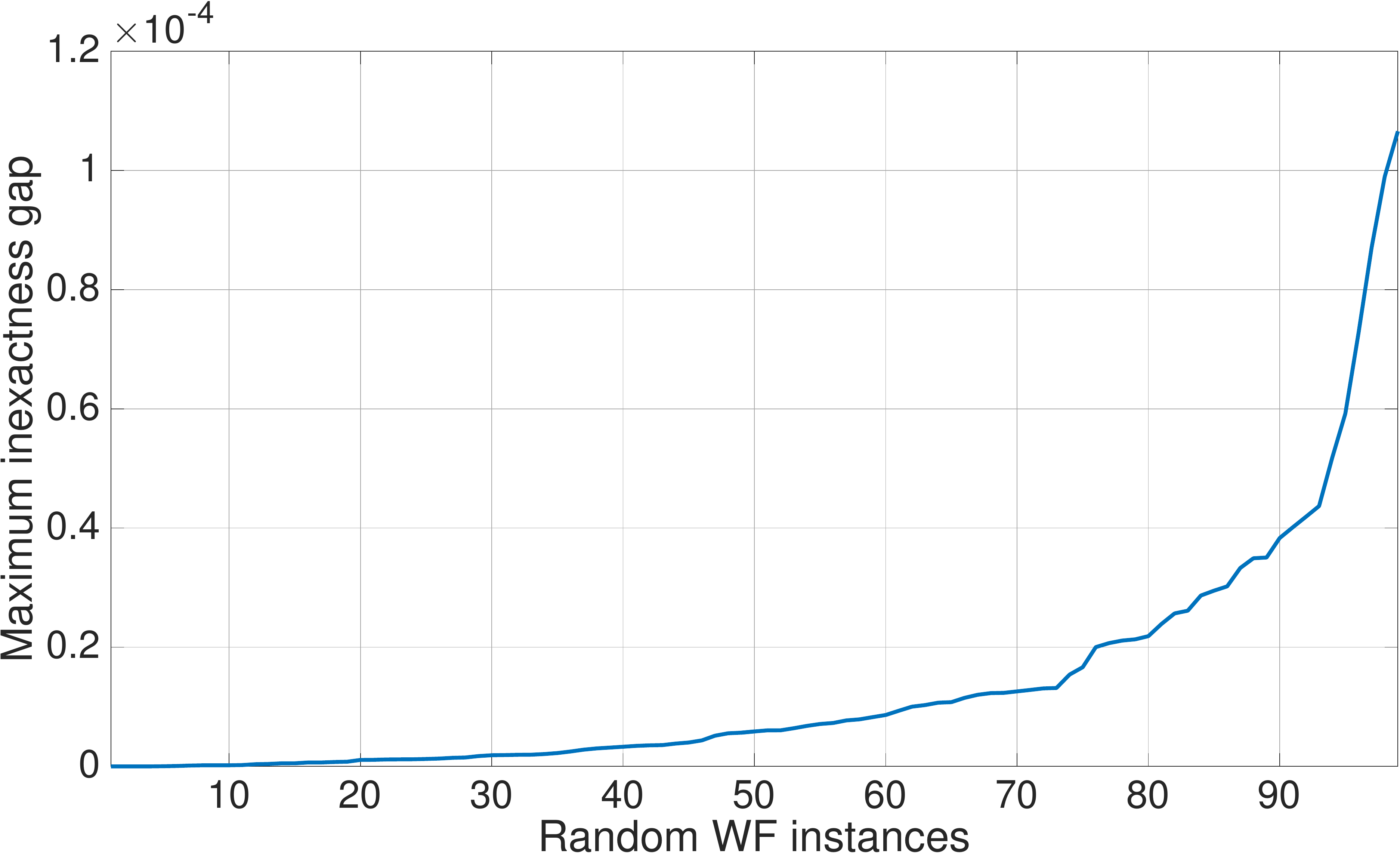}
	%\hspace*{-3em}
	\includegraphics[scale=0.21]{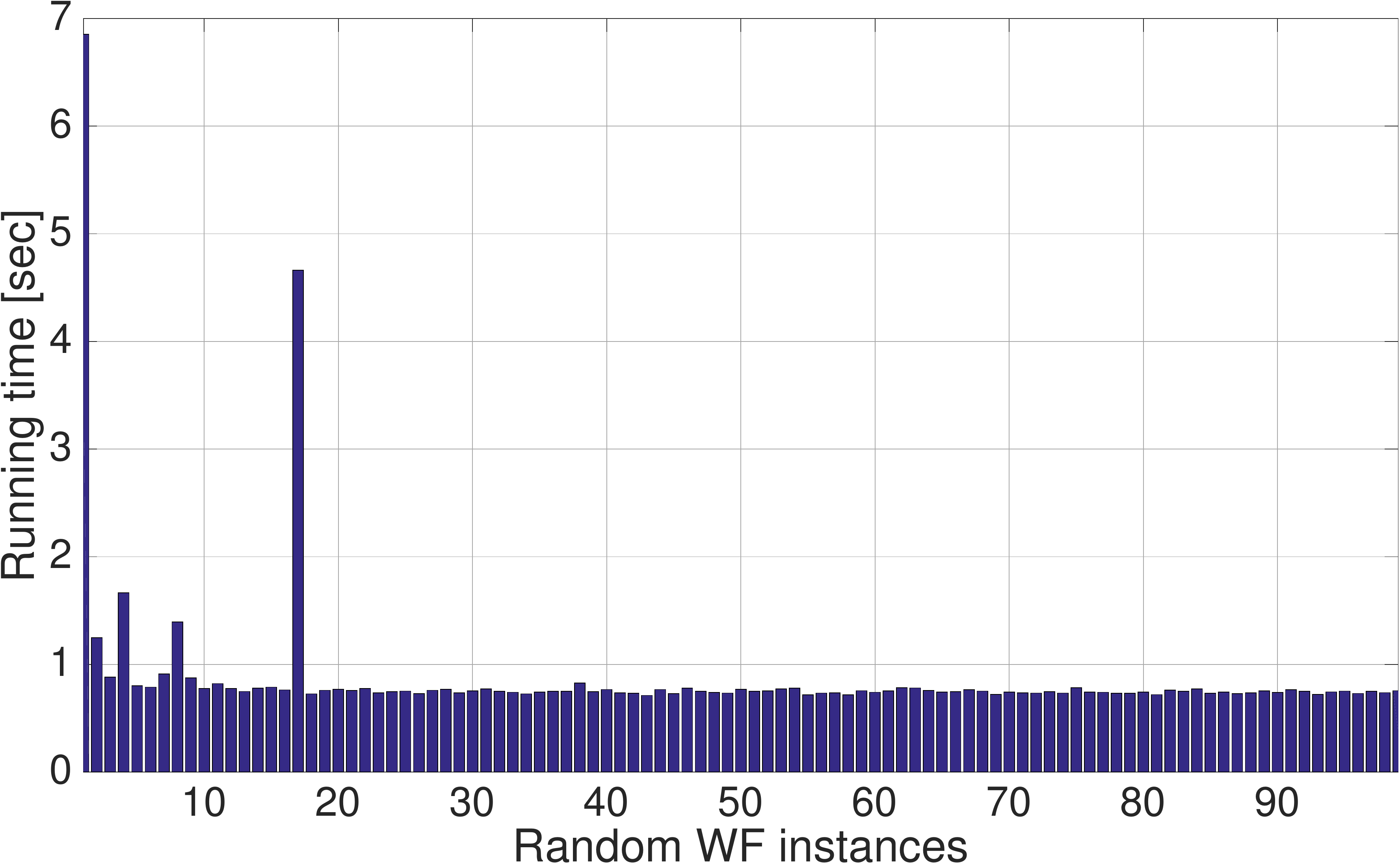}
	%\vspace*{-1em}
	\caption{\emph{Top:} Maximum inexactness gap for random WF instances with all pumps of Fig.~\ref{fig:WDS20} running; \emph{Bottom}: Running times for random WF instances.}
	\label{fig:new}
\end{figure}

To evaluate whether the status of the pumps has any significant effect on runtime, we conducted the previous tests on the WDS of Fig.~\ref{fig:WDS20} with all the pumps running. The inexactness gap and running times were found similar to the tests with pump statuses randomly chosen. Interestingly, the inexactness gap and running times for the new tests did not show significant differences from the tests with random pump statuses; see Fig.~\ref{fig:new}. 99 out of the $100$ random WF instances were solved within the $1$-min deadline with a median running time of $0.75$~sec. Thus, the MI-QCQP \eqref{eq:W2} is in general a powerful WF solver for various network configurations.

%%%%%%%%%%%%%%%%%%%%%%%%%%%%%%%%%%%%%%%%%%%%%%%%%%%%
\section{Conclusions}\label{sec:conclusions}
Using recent tools from graph theory, convex relaxations, energy function-based approaches, and mixed-integer programing, this work has provided a fresh perspective on the physical laws governing water distribution networks. It has been established that the WF problem admits a unique solution even in networks with multiple fixed-pressure nodes and flow-dependent pump models. This WF solution can be provably recovered via a hierarchical stack of WF solvers suitable for different network configurations. Radial networks can be handled by simple convex minimization tasks, whereas networks with cycles call for more elaborate MI-QCQP-based solvers. Nevertheless, numerical tests demonstrate that even the MI-QCQP approach scales well for moderately-sized networks. The MI-QCQP solvers have been derived upon a convex relaxation of the pressure drop equation followed by an objective penalization, an approach that sparked a parallel line of research on the OWF problem~\cite{singh2018optimal}. Network configurations hosting pumps on overlapping cycles remain to be a challenging case.

%%%%%%%%%%%%%%%%%%%%%%%%%%%%%%%%%%%%%%%%%%%%%%%%%%%%
\appendix
\begin{IEEEproof}[Proof of Theorem~\ref{th:wf1}]
The cost of \eqref{eq:W2} can be written as
\begin{eqnarray*}
s(\bh;\bef)=\sum_{(m,n)\in \mcP}\left(h_m-h_n\right)\sign(f_{mn}).
\end{eqnarray*}
To express pressure differences along the flow direction in a compact manner, define the $P\times N$ edge-node incidence matrix $\bA(\bef)$: Its dependence on $\bef$ signifies that the directionality of each edge coincides with the flow directions in $\bef$. Therefore, if the $p$-th row of $\bA(\bef)$ is associated with pipe $p=(m,n)$, then its $(p,k)$ entry is
		\begin{equation*}%\label{eq:Adt}
		A_{p,k}(\bef):=\left\{\begin{array}{ll}
		-\sign^2(f_{mn})+\sign(f_{mn})+1	&,~k=m\\
		\sign^2(f_{mn})-\sign(f_{mn})-1		&,~k=n\\
		0												&,~\textrm{otherwise}.
		\end{array}\right.
		\end{equation*}
For zero flows $(\sign(f_{mn})=0)$, the default pipe direction $(m,n)$ is selected without loss of generality. 
		
Based on $\bA(\bef)$, the pressure differences along the direction of flows can be written as $\bA(\bef)\bh$, and so
\begin{equation}\label{eq:sfun}
s(\bh;\bef)=\boldsymbol{1}^\top\bA(\bef)\bh.
\end{equation}
				
%If \eqref{eq:W1} is feasible, denote its unique solution by $(\bef,\bh)$. Proving by contradiction, suppose $(\tbf,\tbh)$ is a minimizer of \eqref{eq:W2}, which is not feasible for \eqref{eq:W1}. Since both $\bef$ and $\tbf$ satisfy \eqref{eq:wmc} for any incidence matrix and hence for $\bA(\bef)$ too, there exists a nonzero vector $\bn\in\nullspace(\bA(\bef)^\top)$ such that
		%\begin{equation}\label{eq:bn}
		%\tbf=\bef+\bn.
		%\end{equation}
%As shown in \eqref{eq:basis}, the vector $\bn$ can be expressed in terms of a set $\mcL_\mcT$ of fundamental cycles using the flow directions of $\bA(\bef)$. Let us simplify notation as $\bA:=\bA(\bef)$ and $\tbA:=\bA(\tbf)$. Since every edge is assumed to belong to at most one cycle, the set of fundamental cycles is unique irrespective of the spanning tree $\mcT$ selected. Therefore, the set of fundamental cycles in fact contains all cycles in the graph, implying $\mcL_\mcT=\mcL$.

If \eqref{eq:W1} is feasible, denote its unique solution by $(\bar{\bef},\bar{\bh})$. Since the base directionality of the WDS graph $\mcG$ is arbitrary, let it coincide with the water flow directions of $\bar{\bef}$. Using this convention, it follows that $\bar{\bef}\geq\bzero$ and $\bA(\bar{\bef})$ is identical to the base incidence matrix $\bA$. Next, proving by contradiction, suppose $(\tbf,\tbh)$ is a minimizer of \eqref{eq:W2}, which is not feasible for \eqref{eq:W1}. Since both $\bar{\bef}$ and $\tbf$ satisfy \eqref{eq:wmc} for  $\bA=\bA(\bar{\bef})$, there exists a nonzero vector $\bn\in\nullspace(\bA^\top)$ such that
\begin{equation}\label{eq:bn}
\tbf=\bar{\bef}+\bn.
\end{equation}
As shown in \eqref{eq:basis}, vector $\bn$ can be expressed in terms of a set $\mcL_\mcT$ of fundamental cycles using the flow directions of $\bA$. To simplify notation, let us define $\tbA:=\bA(\tbf)$. \color{black}Since every edge is assumed to belong to at most one cycle, the set of fundamental cycles is unique irrespective of the spanning tree $\mcT$ selected. Therefore, the set of fundamental cycles in fact contains all cycles in the graph, implying $\mcL_\mcT=\mcL$.

%Building on \eqref{eq:basis}, consider the $p$-th entry of $\bn$. Because edge $p$ belongs to at most one cycle, either $n_p= 0$ or $n_p=\alpha_\ell n_{\ell,p}$ for a single $\ell\in\mcL$. In the latter case, if $\alpha_\ell<0$, then one can reverse the direction of cycle $\ell$ and substitute $(\alpha_\ell,\bn_\ell)$ in \eqref{eq:basis} with $(-\alpha_\ell,-\bn_\ell)$. Hence, it can be assumed that $\alpha_\ell\geq 0$ for all $\ell\in\mcL$ without loss of generality.

Building on \eqref{eq:basis}, consider a decomposition of $\bn$ as a weighted sum of indicator vectors $\bn_\ell$'s. If there exists a cycle $\ell$, for which $\alpha_\ell<0$, then one can reverse the direction of cycle $\ell$ and substitute $(\alpha_\ell,\bn_\ell)$ in \eqref{eq:basis} with $(-\alpha_\ell,-\bn_\ell)$. Hence, it can be assumed that $\alpha_\ell\geq 0$ for all $\ell\in\mcL$ without loss of generality.
		
Each vector $\bn_\ell$ can be decomposed as $\bn_\ell=\bn_\ell^+-\bn_\ell^-$, where $\bn_\ell^+:=\max\{\bn_\ell,\bzero\}$ and $\bn_\ell^-:=\max\{-\bn_\ell,\bzero\}$. Because every edge belongs to at most one cycle, it holds that 
		\begin{equation}\label{eq:1}
		\boldsymbol{1}=\sum_{\ell \in \mcL}\bn_\ell^++\sum_{\ell \in \mcL}\bn_\ell^-+\bn^0
		\end{equation}
		where the $p$-th entry of vector $\bn^0$ is $1$ if edge $p$ does not belong to any cycle; and $0$, otherwise. Using \eqref{eq:1} in \eqref{eq:sfun}, the objective function of \eqref{eq:W2} becomes
%		\begin{align}\label{eq:sh}
%		s(\bh;\bef)&=\sum_{\ell \in \mcL}(\bn_\ell^+)^\top\bA\bh + \sum_{\ell \in \mcL}(\bn_\ell^-)^\top\bA\bh\nonumber\\
%		&\quad+(\bn^0)^\top\bA\bh.
%		\end{align}
\begin{equation}\label{eq:sh}
s(\bh;\bef)=\sum_{\ell \in \mcL}(\bn_\ell^+)^\top\bA\bh + \sum_{\ell \in \mcL}(\bn_\ell^-)^\top\bA\bh+(\bn^0)^\top\bA\bh.
\end{equation}		
		
Albeit $s(\bh;\bef)$ will be evaluated for different pairs $(\bh;\bef)$, the vectors $\bn_\ell$'s remain unchanged and depend on $\bA$. Based on \eqref{eq:basis} and \eqref{eq:sh}, we will next show that $s(\bar{\bh};\bar{\bef})<s(\tbh;\tbf)$. To do so, we consider the three terms of \eqref{eq:sh} separately.
		
\textit{First summand of \eqref{eq:sh}.} Recall $\bn_\ell^+$ is a binary vector; and the base graph directionality is such that $\bar{\bef}\geq \bzero$. Consider the entries of $\bar{\bef}$ and $\tbf$ for the edges $p$ related to $n_{\ell,p}^+=1$. If $n_{\ell,p}^+=1$, then $\tdf_{p}=\bar{f}_{p}+\alpha_\ell \geq \bar{f}_{p}\geq 0$. In that case, if edge $p=(m,n)$ and relates to a lossy pipe, we get
		\begin{equation}\label{eq:s1pipe}
		(\tdh_{m}-\tdh_{n})\sign(\tdf_{p}) \geq c_p\tdf_{p}^2\geq c_p\bar{f}_{p}^2=(\bar{h}_{m}-\bar{h}_{n})\sign(\bar{f}_{p})
		\end{equation}
		where the first inequality stems from constraint \eqref{eq:relaxed} of \eqref{eq:W2}; the second one from $\tdf_{p}\geq \bar{f}_{p}\geq 0$; and the equality from constraint \eqref{eq:headloss} of \eqref{eq:W1}. 
		
		If edge $p=(m,n)$ relates to a pump, we can also show
		\begin{align}\label{eq:s1pump}
		(\tdh_{m}-\tdh_{n})\sign(\tdf_{p}) & \geq -g_p(\tdf_{p})\geq -g_p(\bar{f}_{p})\nonumber\\
		&=(\bar{h}_{m}-\bar{h}_{n})\sign(\bar{f}_{p}).
		\end{align}
		
Consider cycle $\ell$ and sum up the LHS and RHS of \eqref{eq:s1pipe} or \eqref{eq:s1pump} for all $p$ with $n_{\ell,p}^+=1$ to get
		\begin{equation}\label{eq:term1}
		(\bn_\ell^+)^\top\tbA\tbh\geq(\bn_\ell^+)^\top\bA\bar{\bh},
		\end{equation}
		where the inequality is strict if $\alpha_\ell\neq0$. Summing \eqref{eq:term1} over all cycles provides	 
		\begin{equation}\label{eq:term1b}
		\sum_{\ell \in \mcL}(\bn_\ell^+)^\top\tbA\tbh>\sum_{\ell \in \mcL}(\bn_\ell^+)^\top\bA\bar{\bh},
		\end{equation}
		with strict inequality arising from the fact that not all $\alpha_\ell$'s can be zero for a nonzero $\bn$.

\textit{Second summand of \eqref{eq:sh}.} As explained earlier, if $n_{\ell,p}^+=1$, then $\tdf_{p}=\bar{f}_{p}+\alpha_\ell \geq \bar{f}_{p}\geq 0$ so $\tdf_p$ remains positive. On the other hand, if $n_{\ell,p}^-=1$, then $\tdf_{p}=\bar{f}_{p}-\alpha_\ell< \bar{f}_{p}$. In the latter case, the flow $\tdf_{p}$ has decreased, and its sign may have been reversed to negative. Since all $\bar{f}_p$'s are positive, the flow reversals in $\tdf_p$'s can be modeled as $\bA=\bS \tbA$, where matrix $\bS:=\diag(\sign(\tbf))$ is diagonal with the signs of $\tbf$ on its main diagonal. Since the vectors $\bn_\ell$'s form a basis for $\nullspace(\bA^\top)$, it holds that $\bA^\top\bn_\ell=\bzero$ and so
		\begin{equation}\label{eq:split1}
		(\bn_\ell^+)^\top\bA\bar{\bh}=(\bn_\ell^-)^\top\bA\bar{\bh}.
		\end{equation}
		
		Similar properties hold for $\tbn_\ell:=\bS\bn_\ell$. To see this, the vector $\tbn_\ell$ belongs to $\nullspace(\tbA^\top)$ since
		\[\tbA^\top\tbn_\ell=\tbA^\top\bS\bn_\ell=\bA^\top\bn_\ell=\bzero.\]
		Therefore, we also get that
		\begin{equation}\label{eq:split2}
		(\tbn_\ell^+)^\top\tbA\tbh=(\tbn_\ell^-)^\top\tbA\tbh
		\end{equation}
		where $\tbn_\ell^+:=\max\{\tbn_\ell,\bzero\}$ and $\tbn_\ell^-:=\max\{-\tbn_\ell,\bzero\}$.
		
		By definition of $\tbn_\ell$, if $n_{\ell,p}=1$, then $S_{p,p}=1$ and $\tdn_{\ell,p}=1$. However, if $n_{\ell,p}=-1$, then $S_{p,p}=+1$ or $S_{p,p}=-1$ depending on the sign of $\tdf_p$, and so $\tdn_{\ell,p}=1$ or $\tdn_{\ell,p}=-1$. It therefore follows that
		\begin{subequations}\label{seq:ordering}
			\begin{align}
			\bn_\ell^+ &\leq \tbn_\ell^+\label{seq:ordering:+}\\
			\bn_\ell^- &\geq \tbn_\ell^-.\label{seq:ordering:-}
			\end{align}
		\end{subequations}
		
		Back to the $\ell$-th term of the second summand in \eqref{eq:sh}:
		\begin{align}\label{eq:train}
		(\bn_\ell^-)^\top\tbA\tbh&\stackrel{a}{\geq} (\tbn_\ell^-)^\top\tbA\tbh\stackrel{b}{=}(\tbn_\ell^+)^\top\tbA\tbh\notag\\
		&\stackrel{c}{\geq} (\bn_\ell^+)^\top\tbA\tbh\stackrel{d}{\geq}(\bn_\ell^+)^\top\bA\bar{\bh}\stackrel{e}{=}(\bn_\ell^-)^\top\bA\bar{\bh}.
		\end{align}
%		\begin{align}\label{eq:train}
%		(\bn_\ell^-)^\top\tbA\tbh&\stackrel{a}{\geq} (\tbn_\ell^-)^\top\tbA\tbh
%		\stackrel{b}{=}(\tbn_\ell^+)^\top\tbA\tbh
%		\stackrel{c}{\geq} (\bn_\ell^+)^\top\tbA\tbh\notag\\
%		&\stackrel{d}{>}(\bn_\ell^+)^\top\bA\bh
%		\stackrel{e}{=}(\bn_\ell^-)^\top\bA\bh.
%		\end{align}
		where $(a)$ stems from \eqref{seq:ordering:-}; $(b)$ from \eqref{eq:split2}; $(c)$ from \eqref{seq:ordering:+}; $(d)$ from \eqref{eq:term1}; and $(e)$ from \eqref{eq:split1}. Summing \eqref{eq:train} over $\ell\in\mcL$
		\begin{equation}\label{eq:term2}
		\sum_{\ell \in \mcL}(\bn_\ell^-)^\top\tbA\tbh>\sum_{\ell \in \mcL}(\bn_\ell^-)^\top\bA\bar{\bh},
		\end{equation}
		where the strict inequality stems from that argument $(d)$ in \eqref{eq:train} is strict for $\alpha_\ell\neq0$, and not all $\alpha_\ell$'s can be zero.
		
		\textit{Third summand of \eqref{eq:sh}.} The third summand sums up the pressure differences along the direction of flow for all edges not lying in any cycle. If edge $p=(m,n)$ belongs to this case (i.e., $n_p^0=1$), then $\bar{f}_p=\tdf_p$ and so as in \eqref{eq:s1pipe} we get
		\begin{equation*}
		(\tdh_{m}-\tdh_{n})\sign(\tdf_{p})\geq c_p\tdf_{p}^2= c_p \bar{f}_{p}^2=(\bar{h}_{m}-\bar{h}_{n})\sign(\bar{f}_{p}).
		\end{equation*}
		Summing up over all edges with $n_p^0=1$ yields
		\begin{equation}\label{eq:term3}
		(\bn^0)^\top\tbA\tbh\geq (\bn^0)^\top\bA\bar{\bh}.
		\end{equation}
		
Adding \eqref{eq:term1}, \eqref{eq:term2}, and \eqref{eq:term3} by parts gives $s(\tbh;\tbf)>s(\bar{\bh};\bar{\bef})$. This contradicts that $(\tbf,\tbh)$ is a minimizer of \eqref{eq:W2} since $(\bar{\bef},\bar{\bh})$ is feasible for \eqref{eq:W2}, and concludes the proof.
\end{IEEEproof}

%%%%%%%%%%%%%%%%%%%%%%%%%%%%%%%%%%%%%%%%%%%%%%%%%%%%
\balance
\bibliography{myabrv,water,gas}
\bibliographystyle{IEEEtran}

\begin{IEEEbiography}[{\includegraphics[width=1in,height=1.25in,clip,keepaspectratio]{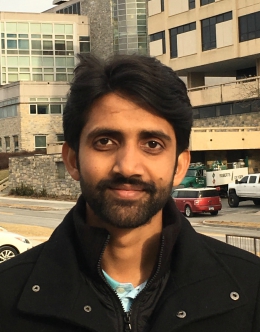}}] {Manish K. Singh} received the B.Tech. degree from the Indian Institute of Technology (BHU), Varanasi, India, in 2013;  and the M.S. degree from Virginia Tech, Blacksburg, VA, USA, in 2018; both in electrical engineering. During 2013-2016, he worked as an Engineer in the Smart Grid Dept. of POWERGRID, the central transmission utility of India. He is currently pursuing a Ph.D. degree at Virginia Tech. His research interests are focused on the application of optimization, control, and graph-theoretic techniques to develop algorithmic solutions for operation and analysis of water, natural gas, and electric power systems.
\end{IEEEbiography} 

\begin{IEEEbiography}[{\includegraphics[width=1in,height=1.25in,clip,keepaspectratio]{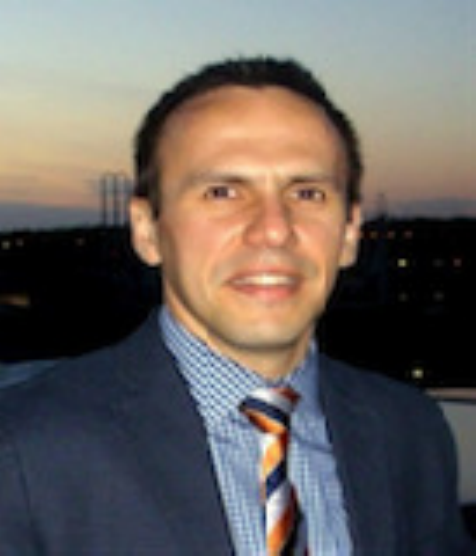}}] {Vassilis Kekatos} (SM'16) is an Assistant Professor with the Bradley Dept. of ECE at Virginia Tech. He obtained his Diploma, M.Sc., and Ph.D. from the Univ. of Patras, Greece, in 2001, 2003, and 2007, respectively, all in computer science and engineering. He is a recipient of the NSF Career Award in 2018 and the Marie Curie Fellowship. He has been a research associate with the ECE Dept. at the Univ. of Minnesota, Minneapolis. During 2014, he stayed with the Univ. of Texas at Austin and the Ohio State Univ. in Columbus as a visiting researcher. His research focus is on optimization and learning for future energy systems. He is currently serving in the editorial board of the IEEE Trans. on Smart Grid.
\end{IEEEbiography}

\end{document}